\documentclass[11pt]{amsart}

\newtheorem{thm}{Theorem}[section]   
\newtheorem{cor}[thm]{Corollary}     
\newtheorem{lem}[thm]{Lemma}         
\newtheorem{prop}[thm]{Proposition}  

\theoremstyle{definition}
\newtheorem{defn}[thm]{Definition}   

\theoremstyle{remark}
\newtheorem{rem}[thm]{Remark}        
\newtheorem{ex}[thm]{Example}        

\numberwithin{equation}{section}     

\newcommand{\thmref}[1]{Theorem~\textup{\ref{#1}}}
\newcommand{\corref}[1]{Corollary~\textup{\ref{#1}}}
\newcommand{\lemref}[1]{Lemma~\textup{\ref{#1}}}

\renewcommand{\H}{\mathcal H}

\newcommand{\B}{\mathcal B}

\renewcommand{\L}{\mathcal L}
\renewcommand{\O}{\mathcal O}
\newcommand{\FF}{\mathcal F}
\newcommand{\F}{\mathbb F}
\newcommand{\A}{\mathcal A}
\newcommand{\D}{\mathcal D}
\newcommand{\C}{\mathbb C}

\newcommand{\Z}{\mathbb Z}

\newcommand{\T}{\mathbb T}
\newcommand{\UM}{UM}
\newcommand{\eps}{\epsilon}

\newcommand{\midtext}[1]{\quad\text{#1}\quad}
\newcommand{\righttext}[1]{\qquad\text{#1 }}

\DeclareMathOperator{\ad}{Ad}
\DeclareMathOperator{\id}{id}

\DeclareMathOperator{\ind}{Ind}
\DeclareMathOperator{\Ind}{Ind}
\DeclareMathOperator{\Aut}{Aut}
\DeclareMathOperator{\Prim}{Prim}
\DeclareMathOperator{\infl}{Inf}

\DeclareMathOperator{\supp}{supp}

\DeclareMathOperator*{\clsp}{\overline{span}}
\DeclareMathOperator*{\spn}{span}

\newcommand{\norm}[1]{\left\| #1 \right\|}
\newcommand{\abs}[1]{\left| #1 \right|}
\newcommand{\imp}[2]{$#1$ -- $#2$}
\newcommand{\case}[1]{\text{if $#1$ (and $0$ else)}}
\newcommand{\inner}[1]{\left\langle #1 \right\rangle}
\newcommand{\lip}[2]{
  {\vphantom\langle}_{#2}\!\! \inner{#1}}
\newcommand{\rip}[2]{
  \inner{#1}_{\!{#2}}  }

\newcommand{\pb}{q^*}
\newcommand{\what}{\widehat}

\begin{document}

\title[Induced coactions]{Induced coactions of discrete groups on
$C^*$-algebras}

\author{Siegfried Echterhoff}
\address{Fachbereich 17, University of Paderborn,
33095 Paderborn, Germany}
\email{echter@math.uni-paderborn.de}

\author{John Quigg}
\address{Department of Mathematics\\Arizona State University\\
Tempe, Arizona 85287}
\email{quigg@math.la.asu.edu}

\thanks{This research is partially supported by National Science
Foundation Grant No. DMS9401253}

\subjclass{Primary 46L55}

\date{December 30, 1997}

\begin{abstract}
Using the close relationship between coactions 
of discrete groups and Fell bundles, we introduce
a procedure for inducing a $C^*$-coaction 
$\delta:D\to D\otimes C^*(G/N)$ of a quotient group
$G/N$ of a discrete group $G$ to a $C^*$-coaction 
$\Ind\delta:\Ind D\to \Ind D\otimes C^*(G)$ of $G$.
We show that induced coactions behave in many respects
similarly to induced actions. In particular, as an analogue of
the well known 
imprimitivity theorem for induced actions
we prove that the crossed products
$\Ind D\times_{\Ind\delta}G$ and $D\times_{\delta}G/N$
are always Morita equivalent. 
We also obtain  nonabelian
analogues of a theorem of Olesen and Pedersen which 
show that there is a duality between
induced coactions and twisted actions in
the sense of Green. We further investigate amenability of Fell bundles
corresponding to induced coactions.
\end{abstract}

\maketitle

\section{Introduction}

One of the most important constructions in ergodic
theory and dynamical systems is the construction
of an induced  action (or induced flow):
if $H$ is a closed subgroup of the
group $G$ and $Y$ is an
$H$-space, then the induced $G$-space
$G\times_HY$ is defined as the quotient of $G\times Y$
with respect to the equivalence relation
$(s,hy)\sim (sh,y)$ for $h\in H$, and $G$-action
given by translation on the first factor.
The induced  $G$-action behaves
in almost all respects similarly to
the original $H$-action on $Y$,
and the theory is particularly useful when it is
possible to identify
a given $G$-space as one which is induced from
a more manageable $H$-space.

The analogue of the induced $G$-space in the
theory of $C^*$-dynamical systems is the induced
$C^*$-algebra $\Ind_H^GA$ together
with the induced action $\Ind\alpha$, where
we start with an action $\alpha:H\to \Aut(A)$.
Needless to  
say,
if $A=C_0(Y)$, then
$\Ind_H^GA=C_0(G\times_HY)$.
As for $G$-spaces, the importance of this
construction comes from the fact that
induced actions enjoy in most respects the same
properties as the original ones.
The most important manifestation of this statement
is certainly Green's imprimitivity theorem
(see \cite[Theorem 17]{gre:local}),
which implies that
the crossed product $\Ind_H^G A\times_{\Ind\alpha}G$ of the
induced system is always Morita equivalent to the crossed
product $A\times_{\alpha}H$ of the original system.
To see the importance of this result,  note that Morita
equivalent algebras have naturally homeomorphic
representation spaces and the same $K$-theory.

In this paper we are concerned with the question whether a similar
theory of induced algebras can be obtained in the
theory of coactions of locally compact groups.
Recall that the theory of coactions of a group $G$
(or rather of the group $C^*$-algebra $C^*(G)$
equipped with a natural comultiplication)
is in a natural way dual to the theory of actions of $G$:
if $\alpha:G\to\Aut(A)$ is an action of $G$ on the $C^*$-algebra $A$,
then there is a canonical coaction
$\widehat{\alpha}$ of $G$ on
$A\times_{\alpha}G$ such that the double crossed product
$A\times_{\alpha}G\times_{\widehat{\alpha}}G$ is stably
isomorphic to $A$ (see \cite{im-tak} and \cite{rae:full}). This generalizes the
Takesaki-Takai  duality theorem for actions of abelian groups,
where $\widehat{\alpha}$ is an action of the dual group
$\widehat{G}$ of $G$.
Conversely, starting with any coaction $\delta$ of $G$ on $A$,
there
exists a dual action $\widehat{\delta}:G\to \Aut(A\times_{\delta}G)$,
and Katayama obtained a similar duality theorem \cite{kat},
which works for all normal coactions (see the preliminary section for the notation).
Of course, in order to develop  the full power of this duality theory,
it is most desirable to have an as complete as possible dual
mirror of the usual constructions for actions. In particular
it would certainly be interesting to have a working notion of induced
coactions and
induced $C^*$-algebras by coactions.
At least if $G$ is discrete we will see here that there is indeed such a theory,
and that it enjoys many properties which are known in the theory
of induced actions.

Our results are based heavily on observations due to the second author,
which connect the theory of coactions of discrete groups to the theory
of Fell bundles (or $C^*$-algebraic bundles) over $G$ \cite{qui:discrete}.
Recall that a Fell bundle $(\A,G)$ over $G$
is a family of Banach spaces $A_s$ for $s\in G$,
together with a multiplication $A_s\times A_t\to A_{st}$ and
an involution $A_s\to A_{s^{-1}}$, which satisfies some further conditions.
The
set $\Gamma_c(\A)$ of sections
of finite support forms a $^*$-algebra,
and a cross sectional algebra $A$ of $(\A,G)$ is a completion
of $\Gamma_c(\A)$ with respect to a given $C^*$-norm.

If $\delta:A\to A\otimes C^*(G)$ is a coaction of a discrete group
$G$ on $A$, then the spectral subspaces $\{A_s:s\in G\}$ 
(i.e., $a_s\in A_s\Leftrightarrow \delta(a_s)=a_s\otimes s$) of 
$\delta$ form a Fell bundle $(\A,G)$ over $G$ and $A$ is
a {\em topologically graded} cross sectional algebra for $\A$
(see \S2 for more details).
Similarly to other situations in the theory of
$C^*$-algebras, there may exist more than one 
$C^*$-norm on $\Gamma_c(\A)$, but if we insist that
the corresponding completions are topologically graded, then
there always exists  a maximal and a minimal one.
We denote the respective cross sectional algebras
by $C^*(\A)$ (for the maximal norm) and $C_r^*(\A)$ 
(for the minimal one; see \cite{exe:amenable} 
for a detailed treatment of this).
Note that both algebras, $C^*(\A)$ and $C_r^*(\A)$, carry  natural
coactions $\delta_{\A}$ and $\delta_{\A}^n$
which are both determined by the property that they map 
$a_s\in A_s$ to the element $a_s\otimes s\in A\otimes C^*(G)$. 
Thus any other coaction $\delta$ lies
``between''  $\delta_{\A}$ and $\delta_{\A}^n$, if $\A$ is the bundle
associated to $\delta$.
Note that just recently, Fell bundles over discrete groups
were studied extensively by several people \cite{qui:discrete, qui-rae:partial,
exe:amenable, ab-exe, exe:partial},
partly due to the discovery that many important $C^*$-algebras
appear as cross sectional algebras of Fell bundles.

Due to the 
above-described 
connection between coactions and Fell bundles
for discrete groups, we are able to use Fell bundles, rather than the coactions
themselves, in order to define induced coactions.
Starting with a Fell bundle $(\D, G/N)$ over a quotient $G/N$  
by a normal subgroup $N$ of a discrete group
$G$,
we define the induced coaction simply as the dual coaction of
the maximal cross sectional algebra $\Ind D:=C^*(q^*\D)$
of the pull-back bundle $(q^*\D, G)$, where
$q:G\to G/N$ denotes the quotient map. 
Note that there is a certain arbitrariness in our definition,
since we could also have taken the dual coaction 
(if it
exists) of
any other topologically graded cross sectional algebra of $q^*\D$
(e.g., $C_r^*(q^*\D)$) instead of the maximal one.
Thus there is no canonical choice unless 
$q^*\D$ is amenable in the sense of Exel \cite{exe:amenable},
which roughly means that all topologically graded cross sectional algebras
are the same. 

We will show that induced coactions behave in almost all respects
similarly
to induced 
actions; for example, if $G$ is abelian, the induced coactions of $G$
correspond exactly to the induced actions of $\widehat{G}$ under the
usual identification between coactions of $G$ and actions of
$\widehat{G}$ (see \S2). 
In \S3 we show that crossed products by coactions of discrete
groups can be realized as cross sectional algebras of certain Fell bundles 
over the transformation groupoid $G\times G$. In fact if
$(\A,G)$ is the Fell bundle associated to the coaction $\delta:A\to A\otimes C^*(G)$,
then $A\times_{\delta}G$ is the 
{\em enveloping $C^*$-algebra} 
of 
$\Gamma_c(\A\times G)$, where $\A\times G$ is the product bundle over 
the groupoid $G\times G$. 
Using 
this result 
we show in \S4
that there is an analogue,  for induced coactions,
of Green's imprimitivity theorem:
if $\delta:D\to D\otimes C^*(G/N)$
is a coaction of $G/N$, then there is a natural Morita equivalence
between the crossed products $D\times_{\delta}G/N$ and
$\Ind D\times_{\Ind\delta}G$. Notice that both crossed products only
depend on the underlying Fell bundles $\D$ and $q^*\D$,
and not on the particular
choices of the cross sectional algebras $D$ and $\Ind D$.

In \S5
 we show that there is an analogue of Olesen
and Pedersen's classical result about twisted group actions 
(see \cite{ole-ped:inner, qui-rae:induce}):  
using  a very useful general characterization of
induced coactions, which is the analogue of the characterization of induced
actions given by the first author in \cite{ech:induce}, we  will see 
that a dual coaction
$\widehat\beta$ of a crossed product
$B\times_{\beta}G$ is induced from a quotient
$G/N$ if and only if  
the action $\beta$
is twisted over $N$ in the sense of Green
\cite{gre:local}.
This leads to a negative result concerning the possibility of a ``Mackey
machine'' for coactions: there is an important feature
of induced actions of compact groups which fails for induced
coactions of discrete groups. Namely, if $\beta$ is an action of a
\emph{compact} 
group $G$ on a $C^*$-algebra $B$
such that $B$ has no proper $G$-invariant ideals,
then $(B,G,\beta)$ is always induced from a system $(A,H,\alpha)$
with $A$ a simple 
$C^*$-algebra; this follows from 
\cite[Theorem]{ech:induce},
since compactness of $G$ guarantees, by
\cite[Lemma 2.1]{ole-ped2}, that $\Prim A$ is equivariantly homeomorphic
to a homogeneous space $G/H$.  For a \emph{discrete} group $G$,
however,  our characterizations of induced coactions allow us to show
that
there exist
numerous examples of $G$-simple coactions which are not induced
(even in the weak sense) from simple coactions!
This drawback of the theory is mainly due to the fact that
the theory of coactions (at least so far) only allows us to look
at quotients by normal subgroups, while for actions we can work
with any closed subgroup of $G$.

Finally, in \S6 we investigate under which conditions
the pull-back bundles $q^*\D$ are amenable. This question is of particular
interest to us, since, as mentioned above,
 only if $q^*\D$ is 
amenable do
we have
a unique choice for our induced algebra $\Ind D$.
In \cite{exe:amenable} Exel introduced a certain approximation property (which
we call property (EP)), which guarantees amenability of
a given Fell bundle $(\A, G)$.
For instance he showed that all Cuntz-Krieger bundles,
which arise from the natural coactions of $\F_n$ on the Cuntz-Krieger
algebras $\O_A$ as found in \cite{qui-rae:partial}, satisfy property (EP),
although the free group $\F_n$ with $n$ generators
is certainly not amenable if $n>1$.
If $\D$ is a bundle over $G/N$, then we will show that
$q^*\D$ satisfies (EP) if $\D$ satisfies (EP)
and $N$ is amenable; the amenability of $N$ is also necessary
for $q^*\D$ 
to satisfy (EP).
Note that as an immediate consequence of
this we see that the Cuntz-Krieger algebras are not induced
from any nontrivial quotient of $\F_n$,
since $\F_n$ does not contain any nontrivial amenable normal subgroup.

This research was conducted while the
second
author visited the University of Paderborn, and he thanks his
hosts Siegfried Echterhoff and Eberhard Kaniuth for their
hospitality.

\section{Preliminaries and basic definitions}
\label{prelim}
Throughout this paper, $G$ will be (except
in certain remarks comparing with other research) a \emph{discrete} group. We are
primarily concerned with coactions of $G$ on $C^*$-algebras, and for these we
adopt the conventions of \cite{qui:discrete} and \cite{qui:fullred}.
We can 
derive
a few benefits from $G$ being discrete: a \emph{coaction}
of $G$ on $A$ is an injective, nondegenerate homomorphism
$\delta\colon A\to A\otimes C^*(G)$ (where here nondegeneracy means
$\clsp{\delta(A)(A\otimes C^*(G))}=A\otimes C^*(G)$)
such that
$(\delta\otimes\id)\circ\delta=(\id\otimes\delta_G)\circ\delta$, where
$\delta_G\colon C^*(G)\to C^*(G)\otimes C^*(G)$ is the homomorphism
defined by $\delta_G(s)=s\otimes s$ for $s\in G$. The \emph{spectral
subspace} of $A$ associated with $s\in G$ is $A_s:=\{a\in
A:\delta(a)=a\otimes s\}$. Since $G$ is discrete, $A$ is the closed span
of the $A_s$. A \emph{covariant representation} of $(A,G,\delta)$ in a
multiplier algebra $M(B)$
is a
pair $(\pi,\mu)$ of nondegenerate homomorphisms of $A$ and $c_0(G)$ 
into $M(B)$ (where, for example, nondegeneracy of $\pi$ means
$\clsp\pi(A)B=B$) such that
\[
\pi(a_s)\mu(\chi_t)=\mu(\chi_{st})\pi(a_s) 
\righttext{for}a_s\in A_s,t\in G, 
\] 
where $\chi_t$ denotes the characteristic function of the 
singleton $\{t\}$. The closed span
$C^*(\pi,\mu):=\clsp\pi(A)\mu(c_0(G))$ is a $C^*$-algebra, and
is called
a \emph{crossed product} for $(A,G,\delta)$ if every covariant
representation $(\rho,\nu)$ factors through $C^*(\pi,\mu)$ in the sense
that there is a homomorphism $\rho\times\nu$ of $C^*(\pi,\mu)$ to
$C^*(\rho,\nu)$ such that $(\rho\times\nu)\circ\pi=\rho$ and
$(\rho\times\nu)\circ\mu=\nu$. All crossed products are isomorphic, and
a generic one is denoted by $A\times_\delta G$, and moreover the
covariant homomorphism generating $A\times_\delta G$ is written
$(j_A,j_G)$. The distinction among the various crossed products is
frequently blurred, and any one of them is referred to as \emph{the}
crossed product. 
The \emph{dual action} of $G$ on the crossed product $A\times_\delta G$
is determined by
$\what\delta_s(j_A(a)j_G(\chi_t))=j_A(a)j_G(\chi_{ts^{-1}})$ for $a\in
A$ and $s,t\in G$.
 
The coaction $\delta$ is called \emph{normal} if $j_A$
is faithful. In any case, there is always a unique ideal $I$ of 
$A$ such that, with $q$
denoting the quotient map from $A$ to $A/I$, the composition
$(q\otimes\id)\circ\delta$ factors through a normal coaction
$\delta^n$, called the \emph{normalization} of $\delta$, on
$A/I$ with the same crossed product as $\delta$, that is, if
$(j_A,j_G^A)$ and $(j_{A/I},j_G^{A/I})$ are the canonical covariant
homomorphisms of $(A,G,\delta)$ into $M(A\times_\delta G)$ and
$M(A/I\times_{\delta^n}G)$, respectively, then $(j_{A/I}\circ q)\times
j_G^{A/I}$ is an isomorphism of $A\times_\delta G$ onto
$A/I\times_{\delta^n}G$. The ideal $I$ coincides with $\ker j_A$, as
well as with $\ker(\id\otimes\lambda)\circ\delta$, where $\lambda$
denotes the left regular representation of $G$.

As shown in \cite{qui:discrete}, for discrete groups coactions are
strongly related to Fell bundles and cross sectional algebras of Fell bundles, for which we adopt the
conventions of \cite{fel-dor} and \cite{exe:amenable}. 
More precisely: if $(A,G,\delta)$
is a coaction, then the spectral subspaces $\{A_s:s\in G\}$ (or, more
properly, the disjoint union of these subspaces) form a Fell bundle
over $G$, which we call the \emph{Fell bundle associated to $\delta$}.

Conversely, if $(\A,G)$ is a Fell bundle, there is a canonical
coaction $\delta_{\A}$, which we will call the \emph{dual coaction},
of $G$ on the \emph{full} cross sectional algebra $C^*(\A)$, determined by 
$\delta_{\A}(a_s)=a_s\otimes s$ for $a_s$ in the fiber $A_s$ of $\A$.
 Throughout this paper, when
we write something like $a_s$ for an element of a Fell bundle, 
we always mean this element is
to be understood to belong to the fiber over $s\in G$.

If
$(\A,G)$ is a Fell bundle over $G$, then a \emph{cross sectional algebra}
$A$ of $(\A,G)$ is simply a completion of $\Gamma_c(\A)$ with respect to any 
given $C^*$-norm. A cross sectional algebra $A$ is called \emph{topologically graded} (see
\cite[Definition 3.4]{exe:amenable}) if there exists a contractive conditional expectation $F:A\to A_e$ which vanishes on each fiber 
$A_s$ for
$s\neq e$ (where we always view the fibers $A_s$ of the bundle
as subspaces of $A$ in the canonical way). 
Exel showed that the full and reduced cross sectional algebras $C^*(\A)$ and
$C_r^*(\A)$ are 
maximal and minimal, respectively, 
among all topologically graded cross sectional algebras of a given
Fell bundle $(\A,G)$. To be more precise: if $A$ is any topologically graded
cross sectional algebra of $\A$, then it follows from 
\cite[Theorem 3.3]{exe:amenable} and the universal property of
$C^*(\A)$ (see \cite[VIII.16.11]{fel-dor}) that the identity map on $\A$
determines 
surjective 
$^*$-homomorphisms 
$$\phi:C^*(\A)\to A,\quad \lambda:A\to C_r^*(\A),\quad\text{and}\quad
\Lambda: C^*(\A)\to C_r^*(A)$$
such that $\lambda\circ \phi=\Lambda$ (the map $\Lambda$ is called the 
\emph{regular representation} of $C^*(\A)$). 
If $\norm{\cdot}_{max}$, $\norm{\cdot}_{\nu}$ and $\norm{\cdot}_{min}$ denote the norms on 
$\Gamma_c(\A)$ coming from viewing $\Gamma_c(\A)$ as a dense subalgebra
of $C^*(\A)$, $A$, and $C_r^*(\A)$, respectively, then the above
result is of course equivalent to saying that 
$\norm{\cdot}_{max}\geq\norm{\cdot}_{\nu}\geq\norm{\cdot}_{min}$.
Thus the topologically graded cross sectional algebras are exactly the completions 
of $\Gamma_c(\A)$ with respect to the $C^*$-norms which lie between 
$\norm{\cdot}_{max}$ and $\norm{\cdot}_{min}$.
Exel calls a Fell bundle $\A$
\emph{amenable} if $C^*(\A)=C^*_r(\A)$ in the sense that the regular
representation of $\A$ is faithful on $C^*(\A)$.  In this case all
topologically graded cross sectional algebras of $\A$ are identical.

If $(A,G,\delta)$ is any coaction, then 
for each $s\in G$ the map
$\delta_s:=(\id\otimes \chi_s)\circ \delta: A\to A$
(where here $\chi_s$, the characteristic function 
of
$\{s\}$,
is regarded as belonging to the
Fourier-Stieltjes algebra $B(G)=C^*(G)^*$, and $\id\otimes\chi_s$ 
is then the slice map of $A\otimes C^*(G)$ into $A$)
is idempotent, with range $A_s$ and kernel containing every $A_t$ for
$t\neq s$. In particular, $\delta_e:A\to A_e$ is a contractive
conditional expectation which vanishes on $A_s$, for all $s\neq e$.
Hence $A$ is a topologically graded cross sectional algebra
of the associated Fell bundle $(\A,G)$.

Now,
\cite[Comment immediately following Definition 3.5]{qui:discrete}
states that the dual coaction $\delta_\A^n\colon a_s\mapsto a_s\otimes
s$ on the \emph{reduced} cross sectional algebra $C^*_r(\A)$ is
(isomorphic to) the normalization of the dual coaction $\delta_\A$ on
$C^*(\A)$. However, there is a subtlety: the constructions of $C^*_r(\A)$
in \cite{qui:discrete} and \cite{exe:amenable} are not quite the
same. So, before we can use the results from both sources, we need to
check that their notions of the reduced $C^*$-algebra of a Fell bundle are
compatible. Namely, we need to know that the kernels in $C^*(\A)$ of the
regular representations of \cite{qui:discrete} and \cite{exe:amenable}
coincide. The conditional expectation $E:=(\delta_\A)_e$ of $C^*(\A)$
onto the
fixed-point algebra $C^*(\A)^{\delta_\A}=A_e$ makes $C^*(\A)$ into
a Hilbert $A_e$-module, as in \cite[Example 6.7]{rie:induced}. Then left
multiplication gives a representation of the $C^*$-algebra $C^*(\A)$
on the Hilbert $A_e$-module $C^*(\A)$, and this in turn gives a
Rieffel inducing map from ideals of $A_e$ to ideals of $C^*(\A)$.
In \cite{qui:discrete} the kernel of the regular representation is
the ideal of $C^*(\A)$ induced from the zero ideal of the fixed-point
algebra $C^*(\A)^{\delta_\A}$. But this coincides with the kernel of the
above representation of $C^*(\A)$ on the Hilbert $A_e$-module $C^*(\A)$,
which is the kernel of the regular representation of \cite{exe:amenable}
(by the proof of \cite[Theorem 3.3]{exe:amenable}). Hence, the
definitions of $C^*_r(\A)$ in \cite{qui:discrete} and
\cite{exe:amenable} are indeed compatible.

The
maps $\phi,\lambda$ and $\Lambda$ considered above
are clearly equivariant with respect to 
the coactions $\delta_{\A}, \delta$ and  $\delta_{\A}^n$
(recall that if $(A,G,\delta)$ and $(B,G,\epsilon)$ are coactions,
then a homomorphism $\phi\colon A\to B$ is called \emph{equivariant}
if $\epsilon\circ\phi=(\phi\otimes\id)\circ\delta$).
Thus we can say that any coaction $\delta:A\to A\otimes C^*(G)$
``lies between'' the dual coaction $\delta_{\A}$ on $C^*(\A)$ 
and its normalization $\delta_{\A}^n$ on $C_r^*(\A)$, 
if $\A$ is the Fell bundle associated to $\delta$.
For reference it is useful to state the following lemma.

\begin{lem}\label{lem-equal}
Let $\delta:A\to A\otimes C^*(G)$ be a coaction of the discrete group 
$G$ and let $(\A,G)$ be the associated Fell bundle.
Let $\delta_{\A}$ and $\delta_{\A}^n$ denote the dual coaction and 
its normalization on $C^*(\A)$ and $C_r^*(\A)$, respectively,
and let $\phi$, $\lambda$ and $\Lambda$ be as above.
Then there are canonical isomorphisms
$$\Ind\phi:C^*(\A)\times_{\delta_{\A}}G\to A\times_{\delta}G,\quad\quad
\Ind\lambda: A\times_{\delta}G\to C_r^*(\A)\times_{\delta_{\A}^n}G,$$
$$\text{and}\quad\Ind\Lambda: C^*(\A)\times_{\delta_{\A}}G
\to C_r^*(\A)\times_{\delta_{\A}^n}G$$
defined by  
\begin{gather*}
\Ind\phi=(j_A\circ \phi)\times j_G^A,\qquad\Ind\lambda= 
(j_{C_r^*(\A)}\circ \lambda)\times j_G^{C_r^*(\A)},\\
\text{and}\qquad
\Ind\Lambda=(j_{C_r^*(\A)}\circ \Lambda)\times j_G^{C_r^*(\A)},
\end{gather*}
respectively.
In particular, $\Ind\Lambda=\Ind\lambda\circ \Ind\phi$ and 
$\delta_{\A}^n$ coincides with the normalization $\delta^n$ of $\delta$.
\end{lem}
\begin{proof}
It follows directly from the equivariance and surjectivity of the maps
$\phi,\lambda$ and $\Lambda$ that the maps $\Ind\phi$, $\Ind\lambda$
and $\Ind\Lambda$ are well defined surjections.  Since
$\Lambda=\lambda\circ \phi$ we also have
$\Ind\Lambda=\Ind\lambda\circ \Ind\phi$.
Since $\delta_{\A}^n$ is the normalization of 
$\delta_{\A}$, it follows from \cite[Corollary 2.7]{qui:fullred}
that $\Ind\Lambda$ is an isomorphism, which then
implies that $\Ind\lambda$ and 
$\Ind\phi$ are also isomorphisms. In particular, it follows that
$\ker j_A=\ker j_{C_r^*(\A)}\circ\lambda=\ker\lambda$ 
(since $j_{C_r^*(\A)}$ is injective by the normality of $\delta_{\A}^n$).
Thus $\delta_{\A}^n$ coincides with the normalization of $\delta$.
\end{proof}

\begin{rem}\label{rem-dual coaction}
In view of the above discussion one could guess that
any topologically graded cross sectional algebra $A$
of a given Fell bundle $(\A,G)$ over the discrete group
$G$ carries a dual coaction $\delta$ which satisfies
$\delta(a_s)=a_s\otimes s$. This is {\em not} the case.

To see a counter example let $G$ be any non-amenable
discrete group such that the direct sum $V=1_G\oplus \lambda_G$
of the trivial representation $1_G$ and the regular representation
$\lambda_G$ of $G$ is not faithful.
Then $V(C^*(G))$ is a topologically graded cross sectional algebra
of the Fell bundle $(\A,G)$ corresponding to $G$ (i.e., $A_s=\C$
for all $s\in G$), since the kernel of $V$ is contained in the kernel
of $\lambda_G$.
Let $U:C^*(G)\to \L(\H)$ be any faithful representation of $G$.
If there were a coaction $\delta$ on $V(C^*(G))$ satisfying
$\delta(a_s)=a_s\otimes s$, this would imply that the unitary
representation $V\otimes U$ of $G$ factors through a faithful
representation of $V(C^*(G))$, i.e., $\ker (V\otimes U)=\ker V$
in $C^*(G)$. But
$V\otimes U=(1_G\oplus\lambda_G)\otimes U=U\oplus (\lambda_G\otimes U)$,
is faithful on $C^*(G)$, while $V$ is not faithful by assumption.

To see that there are numerous examples of groups satisfying the above
property on $1_G\oplus\lambda_G$, let us first note that
any non-amenable group with $1_G\oplus \lambda_G$ faithful
satisfies Kazhdan's
property (T) (i.e., the trivial representation is an isolated point in
$\widehat{G}$).
Since the nonabelian free groups $\F_n$ in $n$ generators
do not satisfy Kazhdan's property (T) (which follows
from the simple fact that
$(\F_n/[\F_n,\F_n])\widehat{\ }=\widehat{\Z^n}=\T^n$),
they all serve as specific examples for our counter example.
Moreover, by a theorem of Fell \cite[Proposition 5.2]{fel:weak}
it is known that for any subgroup $H$ of a discrete group $G$
and any representation $V$ of $H$, $V$ is a direct summand of
$(\ind_H^GV)|_H$, which implies that any
faithful representation of $C^*(G)$ restricts to a faithful representation
of $C^*(H)$.
Thus, if $G$ is a non-amenable group with $1_G\oplus\lambda_G$ faithful,
it follows that $(1_G\oplus\lambda_G)|_H=1_H\oplus \lambda_G|_H$, and hence
$1_H\oplus \lambda_H$ is faithful on $C^*(H)$ for any subgroup
$H$ of $G$, since it follows from
\cite[Addendum of Theorem 1]{herz} that $\ker\lambda_G|_H=\ker\lambda_H$.
In particular, $1_G\oplus \lambda_G$ is not faithful for
any discrete group which contains the free group $\F_2$ as a subgroup.
This shows that any non-amenable group with
$1_G\oplus\lambda_G$ faithful
must indeed be very exotic, and it is certainly an interesting question whether
there exist such groups. 
We are grateful to Alain Valette for some useful comments on this.
\end{rem}

We want to define a notion of induced coactions, dual to the concept of
induced actions. If $N$ is a subgroup of $G$, we can induce an action
of $N$ to an action of $G$,
so dually we should expect to induce a coaction from a quotient
group to the big group. For this we require $N$ to be a normal subgroup.
We don't know yet how to induce coactions in general; for dual coactions
of Fell bundles the way seems fairly clear 
now (given the techniques of the present paper!), 
but for arbitrary
coactions it seems much more difficult. We will develop the theory for
dual coactions of Fell bundles 
over
discrete groups, where the
computations are so much cleaner than for continuous groups. It will be
fairly obvious to the reader that 
some
of what we will do in this paper
can be done for Fell bundles over continuous 
groups, and indeed we plan
to pursue this.

 However, we feel
it is valuable to have the machinery laid out for the case of discrete
groups, since the discrete theory has a flavor all its own. We first define
the Fell bundle which will be associated to the induced coaction. This will
just be the ``Banach
$^*$-algebraic bundle retraction'' by the quotient map $G\to G/N$, as in
\cite[VIII.3.17]{fel-dor}, but we use different notation and terminology:

\begin{defn}
Suppose $(\D,G/N)$ is a Fell bundle over $G/N$, 
where $N$ is a normal subgroup of the discrete
group $G$, and let 
$q\colon G\to G/N$ be the quotient map. We define the \emph{pull-back
Fell bundle over $G$} as
\[ \pb\D=\{(D_{sN},s):s\in G\}. \]
The bundle projection is $(d_{sN},s)\mapsto s$, and we denote the fiber
over $s$ by
$\pb D_s=(D_{sN},s)$. Each fiber $\pb D_s$ is given the Banach
space structure of $D_{sN}$. The multiplication and involution are
defined by
\begin{align*}
(d_{sN},s)(d_{tN},t)&=(d_{sN}d_{tN},st)\\
(d_{sN},s)^*&=(d_{sN}^*,s^{-1}).
\end{align*}
\end{defn}

It is completely routine to verify that the above operations indeed make
$\pb\D$ into a Fell bundle over $G$.

\begin{defn}
Let $(D,G/N,\delta)$ be a coaction, $\D$ the associated Fell
bundle over $G/N$, and $\pb\D$ the pull-back Fell bundle over $G$. We
call the
full cross sectional algebra $C^*(\pb\D)$ the \emph{algebra induced
from $D$} and denote it by $\ind D$, and we call the dual coaction
on $C^*(\pb\D)$ the \emph{coaction induced from $\delta$} and denote it
by $\ind\delta$.
\end{defn}

\begin{rem}\label{rem unique}
Note that the induced algebra and the induced coaction
depend only upon the Fell bundle $\D$; in general $D$ will be some
intermediate algebra between $C^*(\D)$ and $C^*_r(\D)$.
In a sense, our definition of the induced $C^*$-algebra $\Ind D$ above
is somehow artificial: we could have equally well defined
$\Ind D$ as the reduced cross sectional algebra $C^*_r(\pb\D)$,
or any algebra which lies ``between'' the full and the reduced
cross sectional algebras and carries a coaction $\epsilon$ which satisfies
$\epsilon(d_{sN},s)=(d_{sN},s)\otimes s$ for all
$(d_{sN},s)\in \pb D_s$. The only case where there is
really a canonical choice is 
when
$\pb\D$ is amenable in the sense
of Exel \cite{exe:amenable}, since then all cross sectional algebras are
the same. We are going to study this problem in \S6.
Anyway, it follows from Lemma \ref{lem-equal} that the
crossed product $\Ind D\times_{\Ind\delta}G$ is \emph{always}
independent from the choice of the cross sectional algebra for $\pb\D$!
\end{rem}

In view of the above remark 
it makes sense to give also the following

\begin{defn}\label{defn-weak}
Let $(A,G,\delta)$ be a coaction of the discrete group $G$ and 
let $N$ be a 
normal 
subgroup of $G$. We say that
$(A,G,\delta)$ is \emph{weakly induced} from $G/N$ if there
exists a Fell bundle $(\D,G/N)$ such that $(\pb\D,G)$
is isomorphic to the Fell bundle associated to $(A,G,\delta)$.
\end{defn}

Hence a weakly induced coaction is actually induced if and only if
$A$ is equal to the full cross sectional $C^*(\A)$, where
$(\A,G)$ is the Fell bundle associated to $(A,G,\delta)$.

\begin{rem}\label{rem-action}
If $G$ is abelian, then our notion of induced coactions
is the same as the notion of an induced action of the dual
group $\widehat{G}$ of $G$ under the 
one-to-one
correspondence
between coactions of $G$ and actions of $\widehat{G}$.
In order to explain this recall first that if $\alpha:\widehat{G}\to
\Aut(A)$ is an action, then the corresponding coaction $\delta_{\alpha}$ of
$G$ on $A$ is given by
$$\delta_{\alpha}:A\to C(\widehat{G},A);
\bigl(\delta_{\alpha}(a)\bigr)(\chi):=\alpha_{\chi}(a),$$
$\chi\in\widehat{G}$. Here we made the identifications
$C^*(G)\cong C(\widehat{G})$ (via Fourier transform) and
$A\otimes C(\widehat{G})\cong C(\widehat{G},A)$.
Since the Fourier transform of $s\in C^*(G)$ is given by the function
$\chi\mapsto \chi(s)$ on $\widehat{G}$, we see that for $s\in G$
the spectral subspace $A_s$ for
$\delta_{\alpha}$ is given by
$$A_s=\{a\in A: \alpha_{\chi}(a)=
\chi(s)a\;\text{for all}\; \chi\in \widehat{G}\}.$$

Suppose now that $N$ is a subgroup of $G$ and let $\beta:\widehat{G/N}\to
\Aut(B)$
be an action of the subgroup $\widehat{G/N}=N^{\perp}$ of $\widehat{G}$.
The induced
$C^*$-algebra
$\Ind(B,\beta)$ is then defined as
$$\Ind(B,\beta):=\{F\in C(\widehat{G},B):
F(\mu\chi)=\beta_{\bar{\chi}}(F(\mu))\;
\text{for all}\; \chi\in \widehat{G/N}, \mu\in \widehat{G}\},$$
with induced action $\Ind\beta:\widehat{G}\to\Aut(\Ind(B,\beta))$ given by
$$\bigl(\Ind\beta_{\mu}(F)\bigr)(\nu):=F(\bar{\mu}\nu).$$
We claim that the coactions $(\Ind(B,\beta), G,\delta_{\Ind\beta})$
and $(\Ind(B,\delta_{\beta}), G, \Ind\delta_{\beta})$ are isomorphic.
For this it suffices to show that the Fell bundle associated to
$\delta_{\Ind\beta}$ is isomorphic to the pull back of the
bundle $(\B, G/N)$ associated to $\delta_{\beta}$, since by the amenability of
$G$ there is only one topologically graded cross sectional algebra for
this bundle \cite[Theorem 4.7]{exe:amenable}.
Indeed, we claim that the
family of maps
$$\Phi_s:\Ind(B,\beta)_s\to (B_{sN},s); \Phi_s(F_s)=(F_s(1_G), s)$$
is well defined and gives the desired isomorphism of bundles.
To see that it is well defined, it is enough to show that $F_s(1_G)\in B_{sN}$
for all $s\in G$. But
$F_s\in \Ind(B,\beta)_s$ if and only if
$\Ind\beta_{\mu}(F)=\mu(s)F$ for all $\mu\in \widehat{G}$,
so that
$$\beta_{\chi}(F_s(1_G))=F_s(\bar{\chi})=\bigl(\Ind\beta_{\chi}(F_s)\bigr)(1_G)=
\chi(s)F_s(1_G),$$
for all $\chi\in \widehat{G/N}$. Thus $F_s(1_G)\in B_{sN}$.
Each map $\Phi_s$ is norm preserving since
$$\|F_s(\mu)\|=\|\Ind\beta_{\bar{\mu}}F_s(1_G)\|=\|\bar{\mu}(s)F_s(1_G)\|=
\|F_s(1_G)\|.$$
If $b_{sN}\in B_{sN}$ we can define an element
$F_s\in \Ind(B,\beta)_s$ by $F_s(\mu):=\bar{\mu}(s)b_{sN}$,
and then
$\Phi_s(F_s)=(b_{sN},s)$, thus each $\Phi_s$ is surjective.
Finally, since all operations in
$\Ind(B,\beta)$ are pointwise, it follows that the family
of maps $(\Phi_s)_{s\in G}$ respects the bundle operations;
for instance
$$\Phi_{st}(F_sF_t)=(F_s(1_G)F_t(1_G),st)=(F_s(1_G),s)(F_t(1_G),t)=
\Phi_s(F_s)\Phi_t(F_t).$$
\end{rem}

\section{Enveloping algebras and discrete coactions}
\label{env}

In this section we show 
that
for discrete groups
all the $C^*$-calculations we must do are
basically
``pure algebra''.
The following terminology shows what we mean by ``pure algebra''.

\begin{defn}
Let $B_0$ be a $^*$-algebra.
By a \emph{representation} of $B_0$ we always mean a $^*$-homomorphism
of $B_0$ into the bounded operators on a Hilbert space.
We say $B_0$ \emph{has an enveloping
$C^*$-algebra} if the supremum of the $C^*$-seminorms on $B_0$ is
finite, and in this case we call the Hausdorff completion of $B_0$
relative to this largest $C^*$-seminorm the \emph{enveloping
$C^*$-algebra} of $B_0$.
\end{defn}

Thus, if $B_0$ is a $^*$-subalgebra of a $C^*$-algebra $B$, and if every
representation of
$B_0$ is bounded in the norm inherited from $B$, then the closure of
$B_0$ in $B$ is the enveloping
$C^*$-algebra of $B_0$.
For example, we have the well known

\begin{lem}
If $(\A,G)$ is a Fell bundle, then
$C^*(\A)$ is the enveloping $C^*$-algebra of
$\Gamma_c(\A)$.
\end{lem}

\begin{proof}
This follows immediately from the observation that a representation
of $\Gamma_c(\A)$ will be pointwise continuous into the operator
norm topology, hence
continuous for the inductive limit topology.
\end{proof}

The notion of enveloping algebras will be very convenient for
crossed products by coactions: let
$(A,G,\delta)$ be a coaction, and let $\A$ be the associated
Fell bundle. Let $\A\times G$ be the product \emph{Banach} bundle over
$G\times G$. Then $\A\times G$ embeds in the crossed product
$A\times_\delta G$ by identifying $(a_s,t)$ with $j_A(a_s)j_G(\chi_t)$,
and the algebraic operations become
\begin{align}
\label{gpd-bund}
\begin{split}
(a_s,t)(a_u,v)&=(a_sa_u,v)\qquad\case{t=uv}\\
(a_s,t)^*&=(a_s^*,st).
\end{split}
\end{align}
Thus, $\A\times G$ acquires the structure of a Fell bundle \emph{over the
groupoid $G\times G$} in the sense of \cite{kum:fell}, where $G\times G$
is given the transformation group groupoid structure associated with
the action of $G$ on itself by left translation:
\[ (s,tr)(t,r)=(s,r)\midtext{and}(s,t)^{-1}=(s^{-1},st). \]
Although the theory of Fell bundles over groupoids is still in its
infancy, it seems reasonable to expect that many of the properties of
Fell bundles over groups will carry over. In particular, the
$^*$-algebra of finitely supported sections should have 
an enveloping $C^*$-algebra. 
We only need this for the above special 
case (see \corref{cor-bdd1} below).

We will
need to know that
\begin{equation}\label{eq-norm} 
\norm{(a_s,t)}:=\norm{a_s}=\norm{j_A(a_s)j_G(\chi_t)}
\end{equation}
for all $(a_s,t)\in \A\times G$, i.e., that the above described embedding
of $\A\times G$ into $A\times_{\delta}G$ is isometric.
First of all, if $E$
is any subset of $G$, 
then $\chi_E=\sum_{t\in E}\chi_t$ strictly in
$\ell^\infty(G)=M(c_0(G))$, so
\[
j_G(\chi_E)=\sum_{t\in E}j_G(\chi_t),
\]
strictly in $M(A\times_\delta G)$. Consequently, for $a_s\in A_s$ we have
\[
j_A(a_s)=j_A(a_s)1=j_A(a_s)\sum_{t\in G}j_G(\chi_t)=\sum_{t\in G}j_A(a_s)j_G(\chi_t).
\]
Moreover, since the $j_G(\chi_t)$ are orthogonal projections summing strictly
to $1$ in $M(A\times_\delta G)$, and since $j_A$ is faithful on the unit
fiber algebra $A_e$, for $a_e\in A_e$ we have
\[
\norm{a_e}=\norm{j_A(a_e)}
=\sup_{t\in G}\norm{j_A(a_e)j_G(\chi_t)}.
\]
But $\norm{j_A(a_e)j_G(\chi_t)}$ is independent of $t\in G$, since
$j_A(a_e)j_G(\chi_s) =\what\delta_{s^{-1}t}(j_A(a_e)j_G(\chi_t))$ and 
$\what\delta$ is isometric,
being an automorphism of a $C^*$-algebra. Hence, we get
\[
\norm{a_e}=\norm{j_A(a_e)j_G(\chi_t)}\righttext{for all}a_e\in A_e,t\in G.
\]
Since $a_s^*a_s\in A_e$ for all $a_s\in A_s$,
Equation \eqref{eq-norm} now follows from
\begin{align*}
\norm{(a_s,t)}^2&=\norm{a_s}^2=\norm{a_s^*a_s}=\norm{j_A(a_s^*a_s)j_G(\chi_t)}\\
&=\norm{(j_A(a_s)j_G(\chi_t))^*(j_A(a_s)j_G(\chi_t))}\\
&=\norm{j_A(a_s)j_G(\chi_t)}^2.
\end{align*}

We are now going to show that $\Gamma_c(\A\times G)$ 
does have an enveloping $C^*$-algebra,
 namely given by the crossed product
$A\times_{\delta}G$, if
$(\A,G)$ is the bundle associated to a given coaction $\delta$. 
We do it in somewhat more generality, namely for 
(appropriate) dense subspaces of the fibers. Although we do not need
the last whistle in the present paper, we include it since we feel it
will be useful elsewhere.

To prepare for the statement, suppose we have a Fell bundle $(\A,G)$,
and for each $s\in G$ we have a \emph{linear} subspace $B_s$ of $A_s$
such that
\begin{equation}
\label{preFell}
B_sB_t\subset B_{st}\midtext{and}B_s^*\subset B_{s^{-1}}.
\end{equation}
Then $\B:=\bigcup_{s\in G}B_s$ is a subbundle of $\A$ (although not
a Fell bundle, since its fibers may be incomplete). Let $\Gamma_c(\B)$
denote the linear span of $\B$ in $\Gamma_c(\A)$. Then $\Gamma_c(\B)$ is a
$^*$-subalgebra, and conversely any $^*$-subalgebra $B$ of $\Gamma_c(\A)$
which is the linear span of the intersections $B_s:=B\cap A_s$ arises
in this way.

\begin{thm}
\label{dec-bdd}
Let $\A$ be a Fell bundle over the discrete group 
$G$.
Suppose
we have a subbundle $\B$ of $\A$ \textup(with possibly incomplete
fibers\textup), satisfying \eqref{preFell}, such that $A_e$ is the
enveloping $C^*$-algebra of $B_e$ and $C^*(\A)$ is the enveloping
$C^*$-algebra of $\Gamma_c(\B)$.  
Then, regarding $\Gamma_c(\B\times G)$ as
a $^*$-subalgebra of
$\Gamma_c(\A\times G)\subseteq C^*(\A)\times_{\delta_\A}G$ via the inclusion
$\B\hookrightarrow\A$, 
the inductive limit topology on 
$\Gamma_c(\B\times G)$
is stronger than the norm inherited from $C^*(\A)\times_{\delta_\A}G$, and
$C^*(\A)\times_{\delta_\A}G$ is the enveloping $C^*$-algebra of
$\Gamma_c(\B\times G)$.
\end{thm}
\begin{proof}
Note first that 
by Equation \eqref{eq-norm}
each fiber $(B_s,t)$ embeds isometrically into
$C^*(\A)\times_{\delta_\A}G$. 
In particular, this implies that for any finite subset $E$
of $G\times G$ and for any element $b$ of $\Gamma_c(\B\times
G)$ supported in $E$, the norm $\norm{b}$ of $b$ as an element of
$C^*(\A)\times_\delta G$ is bounded by $\max_{(s,t)\in E}\norm{b_{s,t}}\abs{E}$,
where $\abs{E}$ denotes the cardinality of $E$. 
This proves the first part of the theorem.

For the second part, let $\Pi$ be a representation of $\Gamma_c(\B\times
G)$.  We must find a covariant representation $(\pi,\mu)$ of
$(C^*(\A),G,\delta_\A)$ whose integrated form $\pi\times\mu$ is
an extension of $\Pi$. For each $t\in G$ define a representation
$\sigma_{t}$ of the unit fiber algebra $B_e$ of $\B$ by
\[
\sigma_{t}(b_e)=\Pi(b_e,t).
\]
Since $A_e$ is the enveloping $C^*$-algebra of $B_e$, $\sigma_{t}$
extends uniquely to a representation, which we continue to denote by
$\sigma_{t}$, of $A_e$.

We will need to know that $\Pi$ is automatically bounded on each fiber
$(B_s,t)$:
\begin{align*}
\norm{\Pi(b_s,t)}^2
&=\norm{\Pi(b_s,t)^*\Pi(b_s,t)}
=\norm{\Pi(b_s^*,st)\Pi(b_s,t)}
\\&=\norm{\Pi\bigl((b_s^*,st)(b_s,t)\bigr)}
=\norm{\Pi(b_s^*b_s,t)}
\\&=\norm{\sigma_{t}(b_s^*b_s)}
\le\norm{b_s^*b_s}=\norm{b_s}^2=
\norm{(b_s,t)}^2.
\end{align*}
Fix a bounded approximate identity $\{d_i\}$ for $A_e$, and put
\[
p_{t}=\lim\sigma_{t}(d_i),
\]
the limit taken in the weak operator topology. The $p_{t}$'s are
orthogonal projections, and so determine a representation $\mu$ of
$c_0(G)$ such that $\mu(\chi_{t})=p_{t}$. Note that
\begin{align*}
\Pi(b_s,r)p_{t}
&=\Pi(b_s,r)\lim\Pi(d_i,t)
\\&=\lim\Pi(b_s,r)\Pi(d_i,t)
\\&=\lim\Pi(b_sd_i,t)&&\case{r=t}
\\&=\Pi(b_s,r)&&\case{t=r},
\end{align*}
since $\Pi$ is bounded on each fiber $(B_s,r)$.  Also,
\begin{align*}
p_{st}\Pi(b_s,t)
&=\lim\Pi(d_i,st)\Pi(b_s,t)
\\&=\lim\Pi(d_ib_s,t)
\\&=\Pi(b_s,t),
\end{align*}
again by boundedness on the fibers.  Furthermore, one of the above
computations implies $\norm{\Pi(b_s,t)}^2\le\norm{b_s}^2$.  Consequently,
the sum
\[
\pi(b_s):=\sum_{t\in G}\Pi(b_s,t)
\biggl(=\sum_{t\in G}\Pi(b_s,t)p_{t}\biggr)
\]
converges in the weak operator 
topology, 
and
$\norm{\pi(b_s)}\leq\norm{b_s}$ for all $b_s\in B_s$. Standard properties
of the weak operator topology allow one to show
\[
\pi(b_sb_t)=\pi(b_s)\pi(b_t)\midtext{and}
\pi(b_s^*)=\pi(b_s)^*\righttext{for}b_s\in B_s,b_t\in B_t,
\]
so the map $\pi$ extends uniquely to a representation of the $^*$-algebra
$\Gamma_c(\B)$, hence to a representation, which
we still call $\pi$, of $C^*(\A)$, since $C^*(\A)$ is the enveloping
$C^*$-algebra of $\Gamma_c(\B)$. 
The reader
can check that the pair $(\pi,\mu)$ is covariant.  Then we have
\[
(\pi\times\mu)(b_s,t)=\pi(b_s)\mu(\chi_{t})
=\biggl(\sum_{r\in G}\Pi(b_s,r)\biggr)p_{t}
=\Pi(b_s,t),
\]
as desired.
\end{proof}

As a direct consequence of \thmref{dec-bdd} 
together with \lemref{lem-equal} 
we obtain

\begin{cor}\label{cor-bdd1}
Let $(A,G,\delta)$ be a coaction of the discrete group 
$G$ and
$(\A,G)$
the corresponding Fell 
bundle.
Then $A\times_{\delta}G$ is the enveloping
$C^*$-algebra of 
$\Gamma_c(\A\times G)$.
\end{cor}

\section{Imprimitivity theorem}
\label{impthm}
In this section we obtain the anticipated dual mirror of Green's
imprimitivity theorem \cite[Theorem 17]{gre:local}.
We remind the reader that $G$ denotes a \emph{discrete} group and $N$ a
normal subgroup.
Starting 
with a Fell bundle $\D$ over $G/N$, we will construct a
\imp{C^*(\pb\D)\times_{\ind{\delta_\D}}G}{C^*(\D)\times_{\delta_\D} G/N}
imprimitivity bimodule $X$, where ${\delta_\D}$ is the 
dual
coaction.
Moreover, we will make $X$ equivariant for the actions
$\widehat{\ind{\delta_\D}}$ and $\infl\widehat{\delta_\D}$ of $G$, where the latter
action is the inflation to $G$ of the dual action $\widehat{\delta_\D}$ of
$G/N$ on $C^*(\D)\times_{\delta_\D} G/N$.
This is a dual analogue of a result due to Raeburn and the first
author \cite[Theorem 3.3]{ech-rae:induced}, which shows the symmetric
imprimitivity theorem of \cite[Theorem 1.1]{rae:symmetric} (which includes
Green's imprimitivity theorem) is equivariant for suitable coactions.

As usual, we will work with dense subspaces. For
$C^*(\D)\times_{\delta_\D} G/N$ we take the
dense $^*$-subalgebra $C_0:=\Gamma_c(\D\times G/N)$, where we remind
the reader to regard $\D\times G/N$ as a Fell bundle over the groupoid
$G/N\times G/N$, with operations given by \eqref{gpd-bund} (with $G$ replaced by $G/N$).
For $C^*(\pb\D)\times_{\ind{\delta_\D}}G$ we form the corresponding dense
$^*$-subalgebra $B_0:=\Gamma_c(\pb\D\times G)$, except that for
$(d_{sN},s)\in\pb D_s$ and $t\in G$ we
write the pair $((d_{sN},s),t)$ simply as a triple $(d_{sN},s,t)$. For
reference, the operations are
\begin{align*}
(d_{sN},s,uv)(d_{uN},u,v)&=(d_{sN}d_{uN},su,v)\\
(d_{sN},s,t)^*&=(d_{sN}^*,s^{-1},st).
\end{align*}
Also, the actions $\widehat{\ind{\delta_\D}}$ and
$\infl\widehat{\delta_\D}$ 
on the subalgebras $B_0$ and $C_0$ are given
on the generators 
by
\begin{align*}
\widehat{\ind{\delta_\D}}(r)(d_{sN},s,t)&=(d_{sN},s,tr^{-1})\\
\infl\widehat{\delta_\D}(r)(d_{sN},tN)&=(d_{sN},tr^{-1}N).
\end{align*}
Our \imp{B_0}{C_0} pre-imprimitivity bimodule will be
$X_0:=\Gamma_c(\D\times G)$, where $\D\times G$ is the product Banach
bundle over $G/N\times G$. The pre-Hilbert
bimodule operations are defined on the generators by
\begin{align*}
(d_{sN},t)\cdot(d_{uN},vN)&=(d_{sN}d_{uN},t)&&\case{s^{-1}tN=uvN}\\ \rip{(d_{sN},t),(d_{uN},v)}{C_0}
&=(d_{sN}^*d_{uN},u^{-1}vN)&&\case{t=v}\\
(d_{qN},q,r)\cdot(d_{sN},t)&=(d_{qN}d_{sN},qt)&&\case{r=t}\\
\lip{(d_{sN},t),(d_{uN},v)}{B_0}
&=(d_{sN}d_{uN}^*,tv^{-1},v)&&\case{su^{-1}N=tv^{-1}N},
\end{align*}
and then extended bilinearly (or sesquilinearly, as the case may be).
Our action
$\gamma$ of $G$ on $X_0$ is defined on 
the generators by
\[ \gamma_s(d,t)=(d,ts^{-1}). \]
It is easy to check that all operations are continuous in the
inductive limit topologies.

\begin{thm}
\label{imp}
If $(\D,G/N)$ is a Fell bundle, where $N$ is a normal subgroup of the
discrete group $G$, the 
above operations make $X_0$ into a \imp{B_0}{C_0}
pre-imprimitivity bimodule.
Consequently, the completion is a
\imp{C^*(\pb\D)\times_{\ind{\delta_\D}}G}{C^*(\D)\times_{\delta_\D} G/N}
imprimitivity bimodule 
$X$.

Moreover, 
the above formula for $\gamma$ determines an action of $G$
on $X$ which implements a Morita equivalence between the actions
$\widehat{\ind{\delta_\D}}$ and $\infl\widehat{\delta_\D}$.
\end{thm}

\begin{proof}
For the first statement, we must check:
\begin{enumerate}
\item
$X_0$ is a \imp{B_0}{C_0} bimodule;
\item $\lip{b\cdot x,y}{B_0}=b\lip{x,y}{B_0}$
and $\rip{x,y\cdot c}{C_0}=\rip{x,y}{C_0}c$;
\item $\lip{x,y}{B_0}^*=\lip{y,x}{B_0}$
and $\rip{x,y}{C_0}^*=\rip{y,x}{C_0}$;
\item $\lip{x,y}{B_0}$ is linear in $x$
and $\rip{x,y}{C_0}$ is linear in $y$;
\item $x\cdot\rip{y,z}{C_0}=\lip{x,y}{B_0}\cdot z$;
\item $\spn\lip{X_0,X_0}{B_0}$ is dense in $B_0$
and $\spn\rip{X_0,X_0}{C_0}$ is dense in $C_0$;
\item $\lip{x,x}{B_0}\ge 0$ and $\rip{x,x}{C_0}\ge 0$;
\item $\rip{b\cdot x,b\cdot x}{C_0}\le\norm{b}^2\rip{x,x}{C_0}$
and $\lip{x\cdot c,x\cdot c}{B_0}\le\norm{c}^2\rip{x,x}{B_0}$.
\end{enumerate}
The verifications of (i)--(v) are routine, and we will just give one
sample of the computations, leaving the rest to the conscientious
reader.
We will prove (vi)--(vii) in one whack using Rieffel's trick: it
suffices to produce nets in both $B_0$ and $C_0$,
each
term of which is a finite sum of the form $\sum\langle x_i,x_i\rangle$,
which are approximate identities for both the algebras and the module
multiplications, in the inductive limit 
topologies (for example, see the discussion following \cite[Lemma
2]{gre:local}).
We prove (viii) by showing we have homomorphisms of $B_0$ and
$C_0$ into the 
adjointable
operators
on the respective Hilbert modules, which
suffices since by
\corref{cor-bdd1} 
the crossed products mentioned in the statement of
the theorem are the enveloping $C^*$-algebras of $B_0$ and
$C_0$. 
(i) 
follows readily from the
definition. For example, the formula for
the right action of $C_0$ on $X_0$ is obviously bilinear on each product
$(D_{sN},t)\times(D_{uN},vN)$, and then the
extension to $X_0\times C_0$ is bilinear by definition; similarly
for the action of $B_0$. We check associativity of the $B_0$-action,
leaving to
the reader the similar verification for the $C_0$-action:
\begin{align*}\bigl((d_{qN},q,rt)(d_{rN},r,t)\bigr)\cdot(d_{sN},t)
&=(d_{qN}d_{rN},qr,t)\cdot(d_{sN},t)
\\&=(d_{qN}d_{rN}d_{sN},qrt)
\\&=(d_{qN},q,rt)\cdot(d_{rN}d_{sN},rt)
\\&=(d_{qN},q,rt)\cdot\bigl((d_{rN},r,t)\cdot(d_{sN},t)\bigr),
\end{align*}
which shows associativity on the generators, and associativity
follows in general by bilinearity. Similarly, commutativity of the left
and right module multiplications is readily verified on the generators,
hence follows in general by bilinearity. The verifications of
(ii)--(v) are quite similar, and we leave them to the reader.

In order to apply the Rieffel trick, we construct an appropriate
approximate identity for $C_0$,
leaving the (easier) construction for $B_0$ to the reader. Let
$\{d_i\}_{i\in I}$ be a bounded approximate identity for the unit fiber
$D_N$, and let $\FF$ denote the family of finite subsets of $G/N$,
directed by inclusion. For each $F\in\FF$ choose $S_F\subset G$
comprising exactly one element from each coset in $F$. Claim:
\[ \left\{\sum_{t\in S_F}\rip{(d_i^{1/2},t),(d_i^{1/2},t)}{C_0}
\right\}_{(i,F)\in I\times\FF} \]
is an approximate identity for both the algebra $C_0$ and the right module
multiplication of $C_0$ on $X_0$, in the inductive limit topologies. First
of all, note that
\[ \rip{(d_i^{1/2},t),(d_i^{1/2},t)}{C_0}=(d_i,tN). \]
For each generator $(d_{uN},vN)$ we
have
\begin{align*}
\left(\sum_{t\in S_F}(d_i,tN)\right)(d_{uN},vN)
&=\sum_{tN\in F} (d_id_{uN},vN)&&\case{tN=uvN}\\
&=(d_id_{uN},vN)&&\text{whenever}F\ge\{uvN\},
\end{align*}
which tends in the norm of the Banach space $D_{uN}$ to $(d_{uN},vN)$
since $D_ND_{uN}=D_{uN}$. It follows immediately that for all $c\in C_0$
the products $\sum(d_i,tN)c$ converge to $c$ in the inductive limit
topology, and similarly for multiplication on the other side.
For the $C_0$-module multiplication, we similarly have
\[ (d_{sN},r)\cdot\sum(d_i,tN)=(d_{sN}d_i,r)
\righttext{whenever}F\ge\{s^{-1}rN\}, \]
which converges in norm to $(d_{sN},r)$, hence if we replace
$(d_{sN},r)$ by any $x\in X_0$ the convergence will be in the inductive
limit topology, as desired.

We show (viii) for $B_0$, leaving the easier verification for
$C_0$ to the reader. Observe that each generator $(d_{sN},s,t)$ of $B_0$
acts continuously on $X_0$ in the inductive limit topology, hence also
as a bounded operator on the normed space $(X_0,\norm{\cdot}_{C_0})$,
since the inductive limit topology on $X_0$ is stronger than the
pre-Hilbert $C_0$-module norm topology. Then the action of $(d_{sN},s,t)$
extends to the Hilbert module completion $X$, which by the algebraic
properties we saw above gives us a $^*$-homomorphism of $B_0$ into the
$C^*$-algebra of adjointable module maps on $X$. By 
\corref{cor-bdd1}
this homomorphism extends to the enveloping $C^*$-algebra
$C^*(\pb\D)\times_{\ind{\delta_\D}}G$,  hence must be contractive,
and we arrive at the desired inequality.

Finally, $\gamma$ is clearly an algebraic representation of $G$ on
the vector space $X_0$, and straightforward computations (first on the
generators and then extending by bi-additivity) show that
\[
\lip{\gamma_r(x),\gamma_r(y)}{B_0}
=\widehat{\ind{\delta_\D}}(r)\bigl(\lip{x,y}{B_0}\bigr)
\midtext{and}
\gamma_r(x\cdot c)=\gamma_r(x)\cdot\infl\widehat{\delta_\D}(r)(c)
\]
for all $r\in G$, $x,y\in X_0$, and $c\in C$.
Hence, each linear map $\gamma_r$ on $X_0$ extends to $X$, and we get
an algebraic representation
of $G$ on $X$, where the above identities continue to hold.
This is enough to show
$\gamma$ gives a Morita equivalence between the actions
$\widehat{\ind{\delta_\D}}$ and $\infl\widehat{\delta_\D}$.
\end{proof}

\section{Characterizations of induced coactions, the Olesen-Pedersen theorem,
and the Mackey machine}

In this section we want to give a number of useful characterizations of
induced coactions, which are 
analogues
of similar characterizations
in the theory of induced actions. We will then use our results for
a discussion about the possible (or 
impossible)
development of a Mackey machine for coactions.

In \cite{ech:induce} the first author gave a characterization of
induced actions: there must be an equivariant homomorphism of $C_0(G/N)$
into the central multipliers of the algebra carrying the action.
Here we give
a corresponding characterization of induced coactions, involving an
appropriate equivariant homomorphism of $N$. Recall from \cite{fel-dor}
that if $(\A,G)$ is a Fell bundle, a \emph{multiplier} of a fiber $A_s$
(called a ``multiplier of order $s$'' in \cite{fel-dor}) is a pair
$m=(L,R)$ of maps from $\A$ to itself such that $L(A_t)\subset A_{st}$
and $R(A_t)\subset A_{ts}$, and the associativity property
\[ R(a)b=aL(b)\righttext{for}a,b\in\A \]
holds, and one writes $L(a)=ma$ and $R(a)=am$. We write $M(A_s)$ for
the set of all multipliers of $A_s$. The adjoint of $m$ is defined by
$m^*a=(a^*m)^*$ and $am^*=(ma^*)^*$, and $m$ is called \emph{unitary}
if $m^*m=mm^*=1$, the identity element of $M(A_e)$ (and the latter
fortunately agrees with the usual notion of the multipliers of the
$C^*$-algebra $A_e$). We write $UM(A_s)$ for the set of all unitary
multipliers of $A_s$, and $UM(\A)$ for the set of all unitary
multipliers of the bundle $\A$.

\begin{thm}
\label{char}
If $N$ is a normal subgroup of the discrete group $G$, a
Fell bundle $\A$ over $G$ is isomorphic to a pull-back bundle $\pb\D$
for some Fell bundle $\D$ over $G/N$ if and only if there is a
homomorphism $u$ of $N$ into $UM(\A)$ such that
\begin{enumerate}
\item $u_n\in M(A_n)$ for all $n\in N$;
\item $a_su_n=u_{sns^{-1}}a_s$ for all $a_s\in A_s,n\in N$.
\end{enumerate}

Consequently, 
a coaction
$(A,G,\delta)$ is weakly induced from a coaction of
$G/N$ if and only if there is a homomorphism $u$ of $N$ into $UM(A)$
such that $\delta(u_n)=u_n\otimes n$ for $n\in N$ and \textup(ii\textup)
above holds. 
\end{thm}

\begin{proof}
Starting with a Fell bundle $\D$ over $G/N$, it is easy to check that
$u_n:=(1,n)\in M(\pb\D)$ (where the $1$ denotes the identity element of
$M(D_N)$) has the required properties.

Conversely, assume we have a Fell bundle $\A$ and a map $u$ satisfying
(i)--(ii). Let $N$ act on the right of $\A$ by $a\cdot n=au_n$. Our
bundle $\D$ will be the orbit space of this action. Write $[a]$ for the
$N$-orbit of $a\in\A$, and for $s\in G$ define
\[ D_{sN}=\{[a]:a\in A_{sn}\text{ for some }n\in N\}. \]
For each $s$ the orbit map $a\mapsto [a]$ takes $A_s$ bijectively onto
$D_{sN}$, and the resulting Banach space structure on $D_{sN}$ depends
only upon the coset $sN$. More precisely, for $a_s,b_s\in
A_s$, $\lambda\in\C$, and $n\in N$ we have
\begin{align*}
[a_s+b_s]&=[a_su_n+b_su_n]\\
[\lambda a_s]&=[\lambda a_su_n].
\end{align*}
Thus, each fiber $D_{sN}$ has a well defined Banach space structure, and
we get a Banach bundle $\D$ over $G/N$. Define multiplication and
involution in $\D$ by
\begin{align*}
[a][b]&=[ab]\\
[a]^*&=[a^*].
\end{align*}
The reader can easily check that these operations are well defined (for
example, if $a_s\in A_s$, $a_t\in A_t$, and $n,k\in N$, then
\[ \bigl[(a_su_n)(a_tu_k)\bigr]=[a_sa_tu_{t^{-1}ntk}]=[a_sa_t]) \]
and that they give $\D$ a Fell bundle structure (for example, if
$a_s\in A_s$ and
$a_t\in A_t$, then $[a_sa_t]\in D_{stN}$, so
$D_{sN}D_{tN}\subset D_{sNtN}$,
and for the $C^*$-norm property we have
\begin{align*}
\norm{[a_s]^*[a_s]}
&=\norm{[a_s^*][a_s]}
=\norm{[a_s^*a_s]}
=\norm{a_s^*a_s}
\\&=\norm{a_s}^2
=\norm{[a_s]}^2).
\end{align*}
To finish, just check that the map $\phi\colon \A\to\pb\D$ defined by
\[ \phi(a_s)=([a_s],s)\righttext{for}a_s\in A_s \]
is a Fell bundle isomorphism (for example,
\begin{align*}
\phi(a_s)\phi(a_t)
&=([a_s],s)([a_t],t)
=([a_s][a_t],st)
\\&=([a_sa_t],st)
=\phi(a_sa_t)).
\end{align*}The other part follows immediately since a coaction is weakly induced if and
only if the associated Fell bundle is a pull-back.
\end{proof}

\begin{rem}\label{rem-char}
Of course, if  $A=C^*(\A)$, where $\A$ is the bundle
associated to $(A,G,\delta)$ (which is for instance always true
if $\A$ is amenable, which in turn
is always true if $G$ is amenable by \cite[Theorem 4.7]{exe:amenable}),
then we can replace the term ``weakly induced'' by the term ``induced''
in the above theorem.

Note that when $G$ is abelian the above characterization reduces to that of \cite{ech:induce}, since we then have an equivariant nondegenerate
homomorphism of $C_0(\widehat G/N^\perp)$ into the central multipliers of
$C^*(\A)$ (since $u_{sns^{-1}}=u_n$ when $G$ is abelian), and then the
Dauns-Hoffman theorem shows that \cite{ech:induce} applies.
\end{rem}

We now 
give
several applications of the above result
in connection with crossed products by twisted actions in the sense of
Green \cite{gre:local}. 
Let $(B,G,\alpha)$ be an action, and let $(B\times G,G)$ be the associated
semidirect product Fell bundle as in \cite{fel-dor}, so that
$C^*(B\times G)$ is naturally isomorphic to
the crossed product $B\times_\alpha G$ and $C^*_r(B\times G)$ to
$B\times_{\alpha,r}G$.
For reference, the operations on the bundle are
\[ (b,s)(c,t)=(b\alpha_s(c),st)\midtext{and}
(b,s)^*=(\alpha_{s^{-1}}(b)^*,s^{-1}). \]
Recall that a Green-twisted system
$(B,G,N,\alpha,\tau)$ consists of an action $\alpha$ of $G$ on $B$
together with a (strictly continuous) homomorphism $\tau:N\to \UM(B)$
satisfying the equations $\alpha_s(\tau_n)=\tau_{sns^{-1}}$ and
$\alpha_n(b)=\tau_nb\tau_{n^{-1}}$ for all $s\in G, n\in N, b\in B$.

\begin{defn} Let $(B,G,N,\alpha,\tau)$ be a twisted action, and let $N$
act from the right
on the (untwisted) semidirect product bundle $B\times G$ via
\[ (b,s)\cdot n=(b\tau_n,n^{-1}s). \]
The \emph{twisted semidirect product bundle} over $G/N$ is the orbit
space, which we denote by $B\times_N G$.
Let $[b,s]$ denote the $N$-orbit of $(b,s)\in B\times G$. Then the
bundle projection is defined to be $[b,s]\mapsto sN$, and
(since $[b,ns]=[b\tau_n,s]$ for $n\in N$) the fibers are given by
\[ (B\times_N G)_{sN}=\{[b,s]:b\in B\}. \]
\end{defn}

The full and reduced twisted crossed products are naturally isomorphic
to the respective full and reduced cross sectional $C^*$-algebras of the
Fell bundle $B\times_NG$.
Landstad \cite{lan:dual} gave a characterization of reduced crossed
products by actions, in terms of (what we now call) the dual coaction, and it would be useful to have a 
version of this ``Landstad 
duality'' for twisted crossed
products. We prove such a result here when the group is discrete, in
which case it is equivalent to characterize twisted semidirect product
Fell bundles.

\begin{thm}
\label{landstad}
If $N$ is a normal subgroup of the discrete group $G$, a
Fell bundle $\D$ over $G/N$ is isomorphic to a twisted semidirect
product bundle if and only if there is a homomorphism $u$ of $G$ into
$UM(\D)$ with $u_s\in D_{sN}$ 
for all $s\in G$.

Consequently, 
a $C^*$-algebra $D$ is isomorphic to a reduced twisted
crossed product $B\times_{\alpha,\tau,r}G$ by an action of $G$ which is
twisted over $N$ if and only if there are a normal coaction $\delta$
of $G/N$ on $D$ and a homomorphism $u$ of $G$ into $UM(D)$ such that
$\delta(u_s)=u_s\otimes sN$ for $s\in G$.
\end{thm}

\begin{proof}
Given a twisted action $(B,G,N,\alpha,\tau)$, it is easy to check that
taking $u=j_G$, where $j_G:G\to \UM(B\times_{\alpha,\tau,r}G)$
is the canonical homomorphism, gives a map with the
required property.

Conversely, assume we have a Fell bundle $\D$ over $G/N$ and a map $u$
as in the statement of the theorem. Put $B=D_N$. Then $u_sBu_s^*=B$ for
all $s\in G$, so $\alpha_s=\ad u_s$ gives an action $\alpha$ of $G$ on
$B$. Moreover, it is easy to check that the map $\tau:=u|_N$ is a twist
for $\alpha$. So, we have a twisted action $(B,G,N,\alpha,\tau)$. It is
readily verified that the map $\phi\colon B\times_N G\to\D$ defined by
\[ \phi([b,s])=bu_s \]
is well defined (since $\tau=u|_N$) and gives a Fell bundle
isomorphism of the twisted semidirect product bundle $B\times_{N}G$ onto
$\D$, using the equality $D_{sN}=D_Nu_s$ for $s\in G$.

The other part follows immediately, since the coaction $\delta$ is
normal if and only if $D$ is isomorphic to the reduced cross sectional
algebra of the associated Fell bundle.
\end{proof}

\begin{rem}
Alternatively, we could prove the above theorem by pulling back to a
Fell bundle over $G$, suitably carrying along the map $u$, then using
Landstad's original theorem, and finally appealing to a uniqueness
clause in the characterization of $\pb\D$. This would be closer to the
strategy of \cite{qui-rae:induce}, but the above proof is much shorter
and more direct.
\end{rem}

Olesen and Pedersen \cite{ole-ped:inner} proved that if
$(B,G,N,\alpha,\tau)$ is a twisted action of an abelian (not necessarily
discrete) group $G$, then the dual action of
$\widehat G$ on the untwisted crossed product
$B\times_\alpha G$ is induced from the dual action of $N^\perp$ on the
twisted crossed product $B\times_{\alpha,\tau}G$. It often happens
that results about actions of abelian groups can be ``naturally''
viewed as results about coactions, and in this case Raeburn and the
second author proved in \cite{qui-rae:induce} an extension of this
result for nonabelian $G$: if $(B,G,G/N,\delta,\kappa)$ is a twisted
coaction (in the sense of Phillips and Raeburn \cite{phi-rae}), then the
dual action of $G$ on the untwisted crossed product $B\times_\delta G$
is induced from the dual action of $N$ on the twisted crossed product
$B\times_{\delta,\kappa}G$. Here we prove a different sort of (partial,
since our groups are discrete) extension:

\begin{thm}
\label{ole-ped}
If $(B,G,N,\alpha,\tau)$ is a twisted 
action of the discrete group $G$, 
then the semidirect
product bundle $B\times G$ is isomorphic to the pull-back of
the twisted semidirect product bundle $B\times_N G$. Consequently, the
dual coaction of $G$ on the crossed product $B\times_\alpha G$
is induced from the dual coaction of $G/N$ on the twisted crossed
product $B\times_{\alpha,\tau}G$.

Conversely, if $(B,G,\alpha)$ is an action, and if the semidirect
product bundle
$B\times G$ is pulled back from a Fell bundle over $G/N$, or equivalently the
dual coaction $\widehat\alpha$ is induced from a coaction of $G/N$, 
then $\alpha$ is twisted over $N$.
\end{thm}

\begin{proof}
If $(B,G,N,\alpha,\tau)$ is a twisted action,
the reader can easily check that the assignment
\[ (b,s)\mapsto ([b,s],s) \]
gives a Fell bundle isomorphism of $B\times G$ onto $\pb(B\times_N G)$.
The statement concerning induced coactions follows immediately
from the fact that $B\times_{\alpha}G=C^*(B\times G)$.

Conversely, suppose $B\times G$ is isomorphic to a
pull-back. 
Then by \thmref{char}
there is a homomorphism $u\colon
N\to UM(B\times G)$ satisfying (i)--(ii) of that theorem.
The reader
can easily check that the map
\[ \tau_n:=(1,n)u_{n^{-1}} \]
gives a twist for the action $\alpha$.
\end{proof}

A very important step in the modern version of the Mackey machine
for actions is the fact that under favourable circumstances
a $G$-simple action of a group $G$ on a $C^*$-algebra $A$ is automatically
induced from a simple system, i.e., from a system $(D,H,\beta)$
with $D$ simple. This works
especially well for actions of compact groups.
To be more precise: suppose that a
compact group $G$ acts on a
$C^*$-algebra $A$ with Hausdorff primitive ideal space
$\Prim(A)$ such that there exists no nontrivial $G$-invariant ideal of
$A$. Then it follows directly from  the compactness of
$G$ that if we pick any $P\in\Prim(A)$ then $\Prim(A)$ is homeomorphic to
$G/H$ as a $G$-space, where $H=\{s\in G: \alpha_s(P)=P\}$ denotes the
stabilizer of $P$. Applying \cite[Theorem]{ech:induce}, it then follows
that $(A,G,\alpha)$ is isomorphic to the induced system
$(\Ind_H^GD,G,\Ind\beta)$,
where $D=A/P$ and $\beta$ is the associated action (determined by
$\alpha|_H$)
of $H$ on the quotient $A/P$ (note that the Hausdorff assumption
of $\Prim(A)$ could even be 
omitted, by \cite[Lemma 2.1]{ole-ped2}).
The resulting Morita equivalence between $D\times_{\beta}H$ and
$A\times_{\alpha}G$ constitutes one of the main steps
for the Mackey machine for actions (compare with \cite[Theorem 18]{gre:local}).

Since coactions of discrete groups behave somehow similarly to actions
of compact groups (and for abelian $G$ they actually correspond directly
to the actions of the compact group $\widehat{G}$), it would be an easy guess
that a statement similar to the above should be true for coactions
of discrete groups.
But trying to do this we run into trouble very soon:
although one can make sense of a definition of a $G$-simple
coaction (see below), it is certainly not clear at all what quotient of
$G$ should serve as a ``stabilizer'' of a primitive ideal of $A$.
Actually, in what follows next we are going to show that
a result similar to the above can not hold in general for coactions.
We do this by using the previous characterizations of induced coactions.

\begin{defn}[{cf.\ \cite[\S 2]{nil:full},
\cite[Definition 4.1]{lan-phi-rae-sut}}]\label{def-simple}
Let $\delta:A\to A\otimes C^*(G)$
be a coaction of a group $G$ on a $C^*$-algebra $A$.
A closed ideal $I$ of $A$ is
called {\em $\delta$-invariant} if
$\delta(I)(1\otimes C^*(G))=I\otimes C^*(G)$.
$(A,G,\delta)$ is called
{\em $G$-simple} if $A$ has no $\delta$-invariant ideals.
\end{defn}

\begin{rem}\label{rem-invariant}
Of course, the above definition makes sense also
for  coactions of non-discrete groups. Notice that
there are weaker notions of $G$-invariant ideals (and hence,
stronger notions of $G$-simple cosystems).
For instance one could define an ideal of $A$ to be
$G$-invariant if $I=\ker(\phi\otimes \id)\circ \delta$,
where $\phi:A\to A/I$ is the quotient map (e.g., see
\cite[Definition 2.4]{ech-rae:stable}). However, if $G$ is amenable
then it is an easy consequence of \cite[Proposition
4.3]{lan-phi-rae-sut}
that both definitions coincide.
\end{rem}

The following lemma is a special case of \cite[Corollary
3.5]{goo-laz} if $G$ is amenable.

\begin{lem}\label{lem-simple} Let $(B,G,\beta)$ be
a $G$-simple system \textup(i.e., $B$ contains no $G$-invariant closed
ideals\textup).
Then $(B\times_{\beta}G, G,\widehat\beta)$ is
a $G$-simple cosystem.
\end{lem}
\begin{proof}
First note that $\beta$ is $G$-simple if and only if
its 
double dual action $\widehat{\widehat{\beta}}$ is $G$-simple.
This follows easily from the the generalized
duality theorem of Imai and Takai for full crossed products
\cite[Theorem 7]{rae:full}.
Suppose that there exists a nontrivial closed $\widehat\beta$-invariant
ideal $I$ of $A=B\times_{\beta}G$. By
\cite[Propositions 2.1 and 2.2]{nil:full}
there are well defined ``restrictions''
$\widehat\beta_I$ and
$\widehat\beta_{A/I}$ 
of $\widehat\beta$ to $I$ and $A/I$, and
by \cite[Theorem 2.3]{nil:full} we get an exact sequence
$$0\to I\times_{\widehat\beta_I}G\to A\times_{\widehat\beta}G\to
A/I\times_{\widehat\beta_{A/I}}G\to 0$$
with respect to the canonical 
maps, 
which are all
equivariant with respect to the
double dual action
$\widehat{\widehat{\beta}}$ (the last assertion 
following directly 
from the
definition of the maps as given in \cite[Theorem 2.3]{nil:full}).
But this shows that 
$I\times_{\widehat\beta_I}G$ 
is a nontrivial
$\widehat{\widehat{\beta}}$-invariant ideal, which shows
that $\widehat{\widehat{\beta}}$ and hence $\beta$ is not $G$-simple.
\end{proof}

\begin{rem}\label{rem-simple}
For amenable $G$ Gootman and Lazar proved in
\cite[Corollary 3.5]{goo-laz} that
 an action $\beta$ is $G$-simple
if and only if the dual coaction $\widehat{\beta}$ is $G$-simple,
and it follows from \cite[Theorem 3.7]{goo-laz} that the dual
version of this result is also true: if $G$ is amenable,
then a cosystem $(A, G, \delta)$
is $G$-simple if and only if the dual system
$(A\times_{\delta}G,G,\widehat{\delta})$ is $G$-simple.
\end{rem}

As mentioned above, a $G$-simple action of a 
{\em compact} group is always
induced from an action on a simple $C^*$-algebra.
We are now going to show that, unfortunately,
a similar result does \emph{not} hold for coactions of discrete groups.
In fact, we can create a multitude of counterexamples
by using the following easy corollary of Theorem \ref{ole-ped}

\begin{cor}\label{cor-notinduced}
Let $(B,G,\beta)$ be an action with $G$ discrete, and for each $P\in \Prim(B)$ denote by
$S_P\subseteq G$ the stabilizer of $P$ under the corresponding action of
$G$ on $\Prim(B)$. Suppose further that $\bigcap\{S_P: P\in\Prim(B)\}=\{e\}$.
Then 
the
dual coaction $\widehat{\beta}$ is not
\textup(weakly\textup) induced
from any nontrivial quotient of $G$.
\end{cor}
\begin{proof} Assume that there is a nontrivial subgroup $N$ of $G$ such that
$\widehat{\beta}$ is (weakly) induced from $G/N$. Then $\beta$ is twisted over
$N$, 
by Theorem \ref{ole-ped}, which in particular implies that the restriction
of $\beta$ to $N$ is implemented by a  homomorphism $\tau:N\to \UM(B)$.
But this implies that $N$ stabilizes all primitive ideals of $B$,
a contradiction to $\bigcap\{S_P:P\in\Prim(B)\}=\{e\}$.
\end{proof}

\begin{ex}\label{ex-simple}
Let $G$ be a discrete group which has a nontrivial subgroup $H$
such that $\bigcap_{s\in G}sHs^{-1}=\{e\}$. There are many examples
of such groups: e.g., take $G$ any finite simple group
and $H$  any nontrivial subgroup, or take $G=\F_n$ and
$H$ a cyclic subgroup.
Put $B=C_0(G/H)$ and let $\beta:G\to\Aut(B)$ denote the action
given by left translation. Then $\beta$ is a $G$-simple
action, and hence $\widehat{\beta}$ is $G$-simple by Lemma \ref{lem-simple}.
It follows from  Corollary \ref{cor-notinduced} that
$\widehat\beta$ is not induced from any
nontrivial quotient $G/N$, since the stabilizers of 
$\beta$ are just the groups
$sHs^{-1}$, $s\in G$.
However, it follows from Rieffel's version of the
imprimitivity theorem for groups
\cite{rie:induced} 
that $B\times_{\beta}G$ is Morita equivalent to $C^*(H)$,
which in both specific examples mentioned above is not simple
and has Hausdorff primitive ideal space.
\end{ex}

The above discussion shows that if we want to have
a theory for coactions which works similarly to
the full Mackey-Green theory for actions,
we have to introduce something like
``coactions of homogeneous spaces'', i.e., coactions
of quotients $G/H$ by 
not-necessarily-normal
subgroups
$H$ of $G$. Up to now, there is no such theory.

\section{Amenability for pull-back bundles}

Recall that a Fell bundle $(\A, G)$ is called  
amenable
if 
the regular representation
$\Lambda:C^*(\A)\to C_r^*(\A)$ is faithful.
In this section we want to investigate under
what conditions the pull-back bundle $(\pb\D, G)$
of a bundle $(\D, G/N)$ is amenable.
This question is particularly interesting for us, since
only for amenable bundles 
is there 
a unique choice
for the induced algebra $\Ind D$, where $D$ is any cross sectional
algebra of $\D$ which carries a coaction $\epsilon$ of $G/N$ given by
$\epsilon(d_{sN})=d_{sN}\otimes sN$.
Inspired by earlier work of C. Anantharaman-Delaroche, Exel
introduced in \cite{exe:amenable} a certain  approximation property for
Fell bundles,
and he proved that this gives a sufficient condition for $\A$ to be
amenable. Let us recall his condition:

\begin{defn}[{cf. \cite[Definition 4.4]{exe:amenable}}]\label{def-EP}
Let $\A$ be a Fell bundle over
the discrete group
$G$. We say that $\A$ has property (EP), if there exists a net of
functions $f_i:G\to A_e$ with finite supports and satisfying:
\begin{enumerate}
\item
$\sup_{i\in I}\norm{\sum_{s\in G}f_i(s)^*f_i(s)}<\infty$;
\item $\lim_{i\to\infty} \sum_{s\in G} f_i(ts)^*a_tf_i(s)= a_t$
for all $a_t\in A_t$.
\end{enumerate}
\end{defn}

Note that Exel did not require finite supports of the functions
$f_i$ in his definition of the approximation property. But
by \cite[Proposition 4.5]{exe:amenable} his definition is equivalent to
the above. In \cite[Theorem 4.6]{exe:amenable}, Exel showed that property (EP)
implies amenability of $\A$. As the main result of this section we will
prove that a pull-back bundle $(\pb\D,G)$ satisfies (EP) at least if 
$(\D,G/N)$ satisfies (EP) and $N$ is amenable. Moreover, amenability
of $N$ is also a necessary condition. In what follows, if $(\A,G)$ is a 
Fell bundle over $G$ and $H$ is a subgroup of $G$, then $(\A_H, H)$ denotes
the restriction of $\A$ to $H$.
We start with a lemma.

\begin{lem}\label{lem-regular}
Let $\A$ be a Fell bundle over the discrete group $G$ and let $H$ be
a subgroup of $G$. Suppose that $C^*(\A_H)$ is  represented faithfully
as a subalgebra of $B(\H)$ for some Hilbert space $\H$, and let
$\delta:  C^*(\A_H)\to C^*(\A_H)\otimes C^*(H)$ denote the coaction
given by $\delta(h)=a_h\otimes h$. Further,
let $\lambda^G$ denote the left regular representation of $G$ on $l^2(G)$.
Then
$(\id\otimes\lambda^G|_H)\circ\delta$ factors through a faithful representation
of $C_r^*(\A_H)$ on $\H\otimes l^2(G)$.
\end{lem}
\begin{proof}
It follows from 
\cite[Addendum of Theorem 1]{herz} 
that the restriction $\lambda^G|_H$ of
$\lambda^G$ to $H$ is weakly equivalent to the left regular representation $\lambda^H$ of
$H$.
Hence it follows that the kernels of
$\id\otimes \lambda^G|_H$ and $\id\otimes\lambda^H$ coincide
in $C^*(\A_H)\otimes C^*(H)$. Thus we get
$\ker (\id\otimes \lambda^G|_H)\circ \delta
=\ker(\id\otimes\lambda^H)\circ \delta$.
But  $\ker(\id\otimes\lambda^H)\circ \delta$ coincides with 
$\ker j_{C^*(\A_H)}$. Hence, since the normalization of
$\delta$ is the ``dual'' coaction on $C_r^*(\A_H)$,
it follows that $(\id\otimes\lambda^G|_H)\circ \delta$
factors through a faithful map on $C_r^*(\A_H)$ (see \S2).
\end{proof}

In the next proposition we are going to refine Exel's arguments
in order to show that property (EP) for a Fell bundle
$(\A,G)$ even implies amenability of
all  restrictions
$(\A_H, H)$ of $\A$ to subgroups $H$ of $G$.
It is not clear to us whether
amenability, as defined above,  is always inherited 
by
restrictions of Fell bundles
to subgroups --- only if
$\A$ is saturated were we able to show that this is
true for normal subgroups (see Remark \ref{rem-sat} below).
We were also not able (so far) to show that property (EP)
is inherited 
by
restrictions $\A_H$ of $\A$.

\begin{prop}[{cf.\ \cite[Theorem 4.6]{exe:amenable}}]\label{prop-subgroup}
Suppose a Fell bundle $\A$ over the discrete group $G$ 
satisfies
property \textup(EP\textup). Then the restricted bundle $\A_H$ of $\A$
is amenable for any subgroup $H$ of $G$, i.e., the regular representation
$\Lambda^H:C^*(\A_H)\to C^*_r(\A_H)$ is an isomorphism.
\end{prop}
\begin{proof} Let $(f_i)_{i\in I}$ be a net of finitely supported functions
on $G$ satisfying the conditions of Definition \ref{def-EP}.
For each $i\in I$ we define a map
$\Psi_i:\A_H\to\A_H$ by
$$\Psi_i(a_h)=\sum_{s\in G}f_i(hs)^*a_hf_i(s).$$
We claim that each map $\Psi_i$ extends to a bounded linear map
$\Psi_i:C_r^*(\A_H)\to C^*(\A_H)$
satisfying
$\norm{\Psi_i}\leq \norm{\sum_{s\in G} f_i(s)^*f_i(s)}$.
If this is shown, then $(\Psi_i\circ\Lambda^H)(a)$
will converge to $a$ for any section $a\in \Gamma_c(\A_H)$ by property
(2) of Definition \ref{def-EP}, and hence 
for any $a\in C^*(\A_H)$, since the $\Psi_i$
are uniformly bounded by property (1)
of Definition \ref{def-EP}.
Hence, if $a\in \ker \Lambda^H$, then
$a=\lim_{i\to\infty} (\Psi_i\circ\Lambda^H)(a)=0$.

In order to prove the claim assume that
$C^*(\A_H)$ is faithfully represented on
a Hilbert space $\H$ as in Lemma \ref{lem-regular}, and
let $f:G\to A_e$ be any finitely supported map.
Recall also from the lemma that $C_r^*(\A_H)$ is represented faithfully
on $\H\otimes l^2(G)$ via $(\id\otimes\lambda^G|_H)\circ \delta$,
which maps an element $a_h$ in the fibre
$A_h$ of $\A_H$ to the operator $a_h\otimes \lambda^G(h)$.
As in the proof of \cite[Lemma 4.2]{exe:amenable}, we define an operator
$V:\H\to\H\otimes l^2(G)$ by
$$V\xi=\sum_{s\in G} f(s)\xi\otimes\chi_s,$$
where $\chi_s$ denotes the characteristic function of $\{s\}$.
Since
$$\norm{V\xi}^2=\sum_{s\in G}\inner{ f(s)^*f(s)\xi,\xi}\leq
\norm{\sum_{s\in G} f(s)^*f(s)}\norm{\xi}^2,$$
it follows  that $\norm{V}^2\leq \norm{\sum_{s\in G} f(s)^*f(s)}$.
Using the equation $V^*(\xi\otimes\chi_s)=f(s)^*\xi$, we easily compute
\begin{align*}
V^*(a_h&\otimes\lambda^G(h))V\xi
=V^*(a_h\otimes\lambda^G(h))\biggl(\sum_{s\in
G}f(s)\xi\otimes\chi_s\biggr)\\
&=V^*\biggl(\sum_{s\in G}a_hf(s)\otimes\chi_{hs}\biggr)
=\sum_{s\in G}f(hs)^*a_hf(s)\xi.
\end{align*}
Hence, if we define $\Psi:C_r^*(\A_H)\to C^*(\A_H)$ by
$$\Psi(\Lambda^H(a))
=V^*\bigl((\id\otimes\lambda^G|_H)\circ\delta(a)\bigr)V,$$
we see that $\Psi$ is a linear map 
whose norm is
bounded by
$\norm{\sum_{s\in G} f(s)^*f(s)}$, and which maps
$a_h\in A_h$ to $\sum_{s\in G}f(hs)^*a_hf(s)$.
Replacing $f$ by the $f_i$ gives the desired result.
\end{proof}

\begin{rem}\label{rem-sat} Note that if $\A$ is
saturated (i.e., $A_sA_t$ is dense in $A_{st}$ for all $s,t\in G$), 
and $N$ is a normal subgroup of
$G$, then it is not to hard to see that amenability of $\A$ implies amenability
of $\A_N$. The proof involves showing that the restriction of
the regular representation of $\A$ to $\A_N$ has the same kernel as
the regular representation of $\A_N$, and that
$C^*(\A_N)$ imbeds faithfully as a subalgebra of $C^*(\A)$.
The first can be shown by using Exel's description of the
kernel of the regular representation via the canonical conditional
expectation \cite[Proposition 3.6]{exe:amenable}. The second fact follows
from the saturatedness, which allows one to apply
\cite[Corollary XI.12.8]{fel-dor} in order to show that
any faithful representation
of $C^*(\A)$ restricts to a faithful representation of $C^*(\A_N)$.
We do not know whether this result still holds without the
saturatedness condition. Since we don't need these results later,
we omit further details. 

As mentioned earlier, we also do not know
whether property (EP) is inherited to (normal) subgroups.
As for pull backs: we have the strong feeling that 
the pull back $(\pb\D,G)$ of a bundle $(\D, G/N)$ is amenable 
if and only $(\D, D/N)$  and $N$ are, but
all we could prove (so far) is 
\end{rem}

\begin{thm}\label{thm-EP}
Suppose that $G$ is a discrete group with normal subgroup $N$,
and that $\D$ is a Fell bundle over $G/N$.
If $(\D,G/N)$ satisfies Exel's property \textup(EP\textup) and
$N$ is amenable,
then $(\pb\D,G)$ satisfies \textup(EP\textup). Conversely, if
$(\pb\D,G)$ satisfies \textup(EP\textup), then $N$ has to be
amenable.
\end{thm}
\begin{proof}
Suppose first that $\pb\D$ satisfies (EP).
In order to see that $N$ is amenable we consider the
restriction $(\pb\D_N,N)$ of $\pb\D$ to $N$. It follows straight from
the definitions  that $\pb\D_N$ is the trivial bundle
$D_{N}\times N$. Hence, the full cross sectional algebra
is isomorphic to $D_{N}\otimes_{\max}C^*(N)$ and the reduced cross sectional
algebra is isomorphic to $D_{N}\otimes C_r^*(N)$. Further,
the regular representation can be identified with the quotient map
$\id\otimes\lambda^N: D_{N}\otimes_{\max}C^*(N)\to D_{N}\otimes C_r^*(N)$.
Proposition \ref{prop-subgroup} implies that
$D_{N}\times N$ is amenable, from which it follows that
$\id\otimes\lambda^N$ is an isomorphism.
But this implies that
$\lambda^N$ is an isomorphism, too. Hence $N$ is amenable.

Suppose now that $N$ is amenable and $(\D,G/N)$ satisfies (EP).
Let $(f_i)_{i\in
I}$ be a net of finitely supported maps $f_i:G/N\to D_{N}$ which satisfies
the conditions of Definition \ref{def-EP} for $\D$.
Since $N$ is amenable, we can also choose a net $(g_j)_{j\in J}$ of finitely
supported complex-valued functions on $N$ satisfying
$\sum_{n\in N}|g_j(n)|^2\leq 1$ such that the
corresponding matrix coefficients
$$n\mapsto \inner{ \lambda^N(n)g_j,g_j}=\sum_{m\in N}\overline{g_j(nm)}g_j(m)$$
converge pointwise to the trivial function $1_N$ on $N$.
Let $c:G/N\to G$ be any cross section satisfying $c(eN)=e$
(in what follows we will simply write $c(s)$ instead of $c(sN)$).
Every $s\in G$ has a unique decomposition $s=c(s)n_s$ with
$n_s=c(s)^{-1}s\in N$.
For each $(i,j)\in I\times J$ we define
$f_i\times g_j:G\to \pb D_e=(D_{N},e)$ by
$f_i\times g_j(s)=(f_i(sN)g_j(n_s),e)$.

Thus, if $(d_{tN},t)\in (D_{tN},t)$ we compute
\begin{align*}
\sum_{s\in G}f_i\times g_j(ts)^*(d_{tN},t)&f_i\times g_j(s)
=\sum_{sN\in G/N}\sum_{m\in N}f_i\times g_j(tc(s)m)^*(d_{tN},t)f_i\times
g_j(c(s)m)\\
&= \sum_{sN\in G/N}\sum_{m\in N}f_i\times
g_j(c(ts)n_{tc(s)}m)^*(d_{tN},t)f_i\times g_j(c(s)m)\\
&=\biggl(\sum_{sN\in G/N}f_i(tsN)^*d_{tN}f_i(sN)
\sum_{m\in N}\overline{g_j(n_{tc(s)}m)}g_j(m), \,t\biggr).
\end{align*}
Now, if $F$ is any finite subset of $G$ and $\eps>0$ is
given, we may choose $j:=j_{(i,\eps,F)}\in J$ such that
$\abs{1- \sum_{m\in N}\overline{g_{j}(n_{tc(s)}m)}g_{j}(m)}
<\frac{\eps}{c_i},$
for any $sN\in \supp f_i$ and $t\in F$, where we put
$c_i:=\sup_{t\in F}\bigl(\sum_{sN\in G/N} \norm{f_i(tsN)}\norm{f_i(sN)}\bigr)$
(assuming $f_i$ to be non-zero).
Hence, for any $t\in F$ and $(d_{tN},t)\in (D_{tN},t)$ we get
\begin{align*}
&\norm{\sum_{s\in G}f_i\times g_j(ts)^*(d_{tN},t)f_i\times g_j(s)  -(d_{tN},t)
} \\
&=\norm{
\biggl(\sum_{sN\in G/N}f_i(tsN)^*d_{tN}f_i(sN)
\sum_{m\in N}\overline{g_j(n_{tc(s)}m)}g_j(m)- d_{tN},\, t\biggr)}\\
&\leq\norm{
\sum_{sN\in G/N}f_i(tsN)^*d_{tN}f_i(sN)
\sum_{m\in N}\overline{g_j(n_{tc(s)}m)}g_j(m)- d_{tN}}\\
&\leq
\norm{
\sum_{sN\in G/N}f_i(tsN)^*d_{tN}f_i(sN)- d_{tN}
} \\
&\quad\quad\quad\quad+
\norm{
\sum_{sN\in G/N}f_i(tsN)^*d_{tN}f_i(sN)\biggl(1-
\sum_{m\in N}\overline{g_j(n_{tc(s)}m)}g_j(m)\biggr)}\\
&\leq
\norm{
\sum_{sN\in G/N}f_i(tsN)^*d_{tN}f_i(sN)- d_{tN}
}
+
\norm{d_{tN}}\sum_{sN\in G/N}\norm{f_i(tsN)}\norm{f_i(sN)}\frac{\eps}{c_i}\\
&\le\norm{
\sum_{sN\in G/N}f_i(tsN)^*d_{tN}f_i(sN)- d_{tN}
} + \norm{d_{tN}}\eps.
\end{align*}
Now we are ready to produce a net $(h_k)_{k\in K}$ which will enforce (EP) on
$\pb\D$. For this let
$K=I\times (0,\infty)\times \FF$, where $\FF$ denotes the set
of all finite subsets of $G$, equipped with the ordering
$$(i,\eps, F)\geq (i',\eps', F')\Leftrightarrow i\geq i', \eps\leq \eps',\;
\text{and}\; F\supseteq F'.$$
If $k=(i,\eps, F)\in K$ is given, we define
$h_k=f_i\times g_j$, where $j=j_{(i,\eps,F)}$ is chosen as above.
Then it is a direct consequence of the above computations that
$\sum_{s\in G} h_k(ts)^*(d_{tN},t)h_k(s)$ converges to $(d_{tN},t)$
for all $(d_{tN},t)\in \pb\D$; in fact,
if $(d_{tN},t)$ and $\delta>0$ are given,
choose
$k_0= (i_0, \eps_0, F_0)$ such that
$t\in F_0$, $\eps_0< \frac{\delta}{2\norm{d_{tN}}}$ and
$\norm{\sum_{sN\in G/N} f_i(tsN)^*d_{tN}f_i(sN)- d_{tN}
}<\frac{\delta}{2}$ for all $i\geq i_0$.
Then the above computations show that
$$\norm{\sum_{s\in G} h_k(ts)^*(d_{tN},t)h_k(s) - (d_{tN},t)}<\delta$$
for all $k\geq k_0$.
Hence the $(h_k)_{k\in K}$ satisfy condition (i) of Definition \ref{def-EP}.
Finally, it is trivial to check that, if $h_k=f_i\times g_j$, then
\begin{align*}
\norm{\sum_{s\in G}h_k(s)^*h_k(s)}&\leq
\norm{\sum_{sN\in G/N}f_i(sN)^*f_i(sN)}\cdot\biggl(\sum_{n\in N}
|g_j(n)|^2\biggr)\\
&\leq
\norm{\sum_{sN\in G/N}f_i(sN)^*f_i(sN)},
\end{align*}
which proves condition (ii) of Definition \ref{def-EP}.
\end{proof}

\begin{rem}\label{rem-free}
In \cite{qui-rae:partial} Raeburn and the second author discovered
that the Cuntz-Krieger algebras carry a natural coaction of
$\F_n$, $n\geq 2$, and are therefore isomorphic to cross sectional
algebras of the associated Fell bundles $(\A,\F_n)$.
Later \cite{exe:amenable}
Exel showed that all these Fell bundles satisfy property (EP),
and hence are amenable, which shows that the Cuntz-Krieger algebras
are completely determined by the corresponding Fell bundles
(even if they do not satisfy condition (I) as studied by Cuntz and Krieger
in \cite{cun-kri}).
Since any 
nontrivial normal 
subgroup
of $\F_n$ is nonamenable, it follows from Theorem \ref{thm-EP}
that the Cuntz-Krieger algebras are not induced from any
coaction on any proper quotient $G/N$ of $G$.
\end{rem}


\providecommand{\bysame}{\leavevmode\hbox to3em{\hrulefill}\thinspace}

\end{document}